\documentclass[A4,11pt]{article}

\usepackage{color}
\usepackage{array}
\usepackage{ulem}
\usepackage{graphicx}

\usepackage[colorlinks,allcolors=blue]{hyperref}   

\usepackage[ansinew]{inputenc}
\usepackage[english]{babel}
\usepackage{fancyhdr}
\usepackage[active]{srcltx}
\usepackage{amsmath,amssymb}
\usepackage{authblk}
\graphicspath{{../figure/}}

\newcommand{\shcdot}{\!\cdot\!}
\newcommand{\shddot}{\!:\!}

\newcommand{\cA}{{\cal A}}

\newcommand{\aPOD}{^{\tiny\text{POD}}}

\newcommand{\micr}{\,$\mu$m\,}

\newcommand{\veps}{\varepsilon}

\newcommand{\bfm}[1]{\mathbf{#1}}
\newcommand{\bsn}[1]{\boldsymbol{#1}}

\newcommand{\bfa}{\bfm{a}}
\newcommand{\bfb}{\bfm{b}}
\newcommand{\bfc}{\bfm{c}}

\newcommand{\bfu}{\bfm{u}}
\newcommand{\bfw}{\bfm{w}}

\newcommand{\bfS}{\bfm{S}}
\newcommand{\bfX}{\bfm{X}}
\newcommand{\bfV}{\bfm{V}}

\newcommand{\bfP}{\bfm{P}}
\newcommand{\bff}{\bfm{f}}
\newcommand{\bfD}{\bfm{D}}
\newcommand{\bfF}{\bfm{F}}
\newcommand{\bfG}{\bfm{G}}
\newcommand{\bfU}{\bfm{U}}
\newcommand{\tm}{\bfm{1}}    

\newcommand{\bfM}{\bfm{M}}
\newcommand{\bfK}{\bfm{K}}

\newcommand{\bfd}{\bfm{d}}

\newcommand{\bfQ}{\bfm{Q}}
\newcommand{\bfC}{\bfm{C}}
\newcommand{\bfH}{\bfm{H}}

\newcommand{\bfE}{\bfm{E}}

\newcommand{\bfSigma}{\bfm{\Sigma}}

\newcommand{\bfphi}{\bsn{\phi}}

\newcommand{\dS}{{\text d}S}

\newcommand{\dOmega}{{\text d}\Omega}

\usepackage[margin=3cm]{geometry}
\usepackage[font=small,labelfont=bf]{caption}

\begin{document}

\title{Reduced order modeling of nonlinear microstructures through Proper Orthogonal Decomposition}

\author[1]{G.~Gobat}
\author[1]{A.~Opreni}
\author[2]{S.~Fresca}
\author[2]{A.~Manzoni}
\author[1]{A.~Frangi}

\date{}
\affil[1]{Dept.\ of Civil and Environmental Engineering, Politecnico di Milano, P.za Leonardo da Vinci 32, 20133 Milano, Italy}
\affil[2]{MOX - Dept.\ of Mathematics, Politecnico di Milano, P.za Leonardo da Vinci 32, 20133 Milano, Italy}

\maketitle

\begin{abstract}
We apply the Proper Orthogonal Decomposition (POD) method for the efficient simulation of several scenarios  undergone by 
Micro-Electro-Mechanical-Systems, involving nonlinearites of geometric and electrostatic nature. 
The former type of nonlinearity, associated to the large displacements of the devices, 
leads to polynomial terms up to cubic order that are reduced through exact projection onto a low-dimensional subspace spanned by the 
Proper Orthogonal Modes (POMs). 
On the contrary, electrostatic nonlinearities are modeled resorting to precomputed manifolds in terms of the amplitudes of the electrically active POMs. 
We extensively test the reliability of the assumed linear trial space in challenging applications focusing on resonators, micromirrors and arches also displaying internal resonances.
We discuss several options to generate the matrix of snapshots using both classical time marching schemes
and more advanced Harmonic Balance (HB) approaches.
Furthermore, we propose a comparison between the periodic orbits computed with POD 
and the invariant manifold approximated with Direct Parametrization approaches, further 
stressing the reliability of the technique and its remarkable predictive capabilities, 
e.g., in terms of estimation of the frequency response function of selected output quantities of interest.
\end{abstract}

\section{Introduction}
\label{sec:intro}

Although model order reduction methods for structures experiencing large-amplitude vibrations with geometric nonlinearities have been investigated for a long time \cite{rega2005,kerschen2005,gordon2008,mignolet}, they have been only recently applied 
to the analysis of Micro-Electro-Mechanical Systems (MEMS), a class of devices with profound and increasing impact in the consumer and automotive market \cite{ijnm19,jmems20reso}.
MEMS structures are generally actuated near resonance and are subjected to relatively large transformations.  
These effects are strongly enhanced by the fact that MEMS are monolithic devices often packaged in near-vacuum, thus limiting dissipation to negligible levels. As a consequence, they show highly nonlinear dynamical features that are rarely observed at the macro scale, ranging from jump phenomena \cite{nayfeh95}, to bifurcations of solutions \cite{nayfeh89} (e.g.  bistability \cite{krylov2011}), internal resonances and saturation effects \cite{Czaplewski2019,Houri2019,avoort10,ruzziconi2021two}, self-induced parametric amplification \cite{nitzan16} and frequency combs \cite{seshia17,seshia18}.
Furthermore, the nonlinear properties of MEMS can be tailored to yield performance that would not be accessible operating in the linear regime \cite{harvester12}. 
Accurate and predictive modeling needs to account for all these aspects. 

However, relying on Full Order Models (FOMs) for the numerical simulation of the structural behavior of MEMS poses severe computational challenges that have been only partially solved so far. Generally, one is interested primarily in the steady-state periodic response of MEMS as a function of the actuation intensity and frequency, i.e., the so-called Frequency Response Function (FRF) of selected output quantities of interest like, e.g., the maximum midspan deflection of a beam, or the rotation amplitude of a micromirror. Moreover, the actuation can be electrostatic, piezoelectric, or magnetic, according to the considered applications, hence introducing additional sources of nonlinearity. Finally, because of the large quality factors involved, ultimately  leading to long transients, time marching schemes are hardly computationally affordable. 
Recent advancements on the topic have enabled FOMs simulations within reasonable computational times \cite{KerschenHB2015,actuators21}. Geometrical and inertial nonlinearities can be modeled using Harmonic Balance (HB) approaches or shooting techniques, which directly compute the periodic response. However, these approaches entail  huge computational costs when applied to large Finite Element Method (FEM) models of MEMS. This motivates the interest in developing rapid and reliable Reduced Order Models (ROMs) that ensure a fast and accurate estimation of the FRF of structures within time spans that are compatible with industrial design requirements.

A large family of ROMs, which we can refer to as {\it linear} ROMs, gathers Galerkin projections onto low-dimensional linear subspaces. One of the simplest options is to use a selection of linear eigenmodes  and resort to procedures like the STiffness Evaluation Procedure (STEP), first introduced in \cite{rizzi2003} to compute coupling coefficients. However, as recently highlighted in \cite{givois,cyrilstep}, its application to 3D FEM models is critical since it is mandatory to explicitly include all the coupled high frequency (e.g. axial, lateral contraction) linear modes which are usually difficult to identify and costly to compute.  This issue has been overcome by the Implicit Condensation (IC) approach, which has been successfully applied to MEMS only recently \cite{ijnm19,jmems20reso,gyro2021,nicolaidou2020indirect}. In this case, a small subset of linear eigenmodes, known as master modes, is defined to span a  stress manifold that statically condenses all the contributions of high frequency modes.  However, when inertia nonlinearities play a major role or the frequencies of the slave modes are not well separated from the master ones, the method fails.  Another linear ROM relies on {\it Proper Orthogonal Decomposition} (POD) \cite{kerschen2005,alfiobook,lu2019}, which  this contribution focuses on. In this case, basis functions are computed in a data-driven manner, performing the singular value decomposition of a matrix of FOM solutions computed over time, and for suitably sampled parameter values; thanks to SVD, the most relevant contributions to explain the solution variability across  the time span and the parameter space are selected, resorting to an energy measure. Later, a Galerkin projection onto the POD subspace allows us to generate a low-dimensional ROM, which we refer to as POD-Galerkin ROM. In this contribution, we show how this approach can overcome the limitations shown by other linear methods.

An alternative to linear ROMs is provided by a different class of methods, which we can refer to as {\it nonlinear ROMs}. Among them, a further classification can be proposed with respect to the assumptions used in their derivation. Modal Derivatives and the related Quadratic Manifold approach \cite{md2016,QM2017,QMgen2017,mahdiabadi2021non} try to define spaces of nonlinear basis functions 
with the key idea of taking into account the amplitude dependence of mode shapes and eigenfrequencies.
However, these functions are assumed to be velocity independent, ultimately introducing model limitations similar to the IC approach.

On the other hand, truly nonlinear reduction methods start by defining a nonlinear relationship between the original coordinates and those of the reduced dynamics, hence providing a more accurate treatment of the nonlinear trajectories and  faster convergence with fewer master modes.  This class of methods resorts to the concept of Nonlinear Normal Mode (NNM),  whose study began with the pioneering work by Rosenberg \cite{rosenberg}. In his work, the NNM was defined as a synchronous vibration of the system. This concept has been later generalized by the notion of invariant manifold \cite{shawNNM,touze2004,touze2006,touze} and spectral submanifold \cite{haller16,haller18}.  While the numerical computation of NNMs for large-scale FEM models has been tackled, e.g., in \cite{KerschenNNM}, the generation of ROMs based on the concept of NNM has been addressed so far for small systems with few degrees of freedom (dofs) and only in very recent contributions \cite{vizza2020,nld21} the technique  has been applied to complex structures involving inertia and geometrical nonlinearities. However, its extension to multiphysics (e.g. electromechanics) has not 
been addressed yet, and poses severe computational challenges.

Early applications of POD \cite{kerschen2002,kerschen2005,sampaio} to elastic structures
with distributed nonlinearities have put in evidence its optimality in the sense that it minimizes the average
distance between the original signal and its reduced linear representation. Indeed, the linear nature of a POD-Galerkin approach can be considered as an advantage since few manipulations are needed to construct the ROM. 
Nevertheless, it  also represents a drawback, because a single, global linear subspace might not be able in principle to describe the nonlinear invariant manifolds \cite{amabili07}
which characterize mechanical structures.
While applications of POD to MEMS \cite{tiso,korvink} have been so far mainly limited 
to linear mechanics, beam theory and optimization problems, in this contribution we focus on the application of POD to highly nonlinear problems, showing its accuracy and computational efficiency.
In particular, different sources of nonlinearities are considered, dealing with large rotations, internal resonances (i.e. nonlinear coupling) and electrostatic forcing. The POD-Galerkin ROM is validated against FOM solutions, and its generalization capabilities over the space of parameters are assessed. In particular, the ROM dynamics is solved by resorting to numerical continuation and bifurcation analysis tools, which give an insight onto the underlying dynamics -- being this latter usually difficult to access in the FOM case because of the high computational cost required. The POD-Galerkin ROM solution is also studied from the perspective of invariant manifold theory. In particular, we compare the periodic orbits obtained from the POD-Galerkin ROM and the invariant manifold approximated with the Direct Parametrization (DP) approach \cite{nld21} applied to the corresponding FEM system, providing a detailed analysis unprecedented for large Finite Element Models.

The structure of the paper is as follows. 
After a short description of the POD-Galerkin framework
in Section~\ref{sec:formulation}, we focus on a series of applications
to MEMS modeled with the FEM, ranging from simple beam resonators to complex micromirrors and 
to arches displaying internal resonances, which put a strain 
on known techniques for geometrical nonlinearities.
In the examples, we discuss the physical rationale behind 
the POD modeling capabilities,  resorting 
to comparisons with the results from Direct Parametrization.
We finally show that the POD-Galerkin ROM  yields great promise also on coupled 
multiphysics applications where few or no alternatives are available.

\section{Formulation}
\label{sec:formulation}
Let us consider the framework of structures subjected to large transformations and small strains. This is the operating range of most microsystems since they are often actuated at resonance and large aspect ratios allow reaching large displacements within the linear elastic range of the material. In this framework, the Saint Venant-Kirchhoff constitutive model \cite{malvern} is the most appropriate choice, and is given by
\begin{equation}
\label{eq:gltensor}
\bfS=\cA:\bfE,
\end{equation}
where $\bfS$ is the second Piola-Kirchhoff strain tensor, $\cA$ the fourth-order elasticity tensor and $\bfE$ the Green-Lagrangian Strain tensor:
\begin{equation}
\bfE =\frac{1}{2}\left(\nabla\bfd+\nabla^T\bfd+\nabla^T\bfd\cdot\nabla\bfd\right);
\end{equation}
here we denote by $\bfd$  the displacement field and by $\nabla(\cdot)$ the (material) gradient defined with respect to the reference configuration. The weak form of the linear momentum conservation law is:
\begin{equation}
\label{eq:PPV}
\begin{split}
\int_{\Omega_0}\!\rho_0\ddot{\bfd}\shcdot\bfw\,\dOmega_0 & +
\int_{\Omega_0}\!\bfP[\bfd]\shddot\nabla^T\bfw\,\dOmega_0  = \\
& \int_{\Omega_0}\!\rho_0\bfF\shcdot\bfw \,\dOmega_0+
\int_{S_T}\!\bff\shcdot\bfw\,\dS_0,  \quad \forall \bfw\in H^1_0(\Omega_0),
\end{split}
\end{equation}
where the integrals are expressed in the reference configuration $\Omega_0$ and $\dot{ \square }$ denotes the time derivative. Here 
$\rho_0$ denotes the initial density,  $\bfP[\bfd] = (\tm+\nabla\bfd)\cdot\bfS$  the 
first Piola-Kirchhoff stress tensor, $\bfF$ the body forces per unit mass, 
$\bff$ the surface tractions prescribed on the surface $S_T$ and $\bfw$  the test velocity selected in $H^1_0(\Omega_0)$, i.e.\ the space of functions with finite energy that vanish on the portion $S_U \subset \partial \Omega_0$ where Dirichlet boundary conditions are prescribed.
In our applications we assume that vanishing
displacements are enforced on $S_U$.
Within the present context, it is worth stressing that eq.\eqref{eq:PPV} exactly accounts for geometric (elastic and inertia) nonlinearities, e.g., large rotations or nonlinear mode coupling.  

The spatial discretization of eq.\eqref{eq:PPV}, e.g.\ by means of finite elements, also including a Rayleigh model damping term, yields to a system of coupled nonlinear differential equations of the following form:
\begin{equation}\label{eq:PPV_d2}
\bfM \ddot{\bfD} + \bfC\dot{\bfD}  + \bfK\bfD + \bfG(\bfD,\bfD) + \bfH(\bfD,\bfD,\bfD) = 
\bfF(\bfD,\bsn{\beta},\omega,t), \qquad t \in (0,T)
\end{equation}
where the vector $\bfD \in \mathbb{R}^n$ collects all unknown displacement nodal values,
$\bfM \in \mathbb{R}^{n \times n}$ is the mass matrix, $\bfC=\omega_0/Q\bfM$ the Rayleigh model mass proportional 
damping matrix -- considering a reference eigenfrequency $\omega_0$ and a quality factor $Q$ -- and  
$\bfF \in \mathbb{R}^{n}$ the nodal force vector which depends on the actuation intensity parameters 
$\bsn{\beta}$, 
the angular frequency of the actuation $\omega$ and in general also on $\bfD$, e.g.\ in electromechanical applications. 
The internal force vector has been exactly decomposed in linear, quadratic, and cubic power terms of the displacement: $\bfK \in \mathbb{R}^{n \times n} $ is the stiffness  matrix related to the linearized system, while $\bfG \in \mathbb{R}^{n}$ and $\bfH \in \mathbb{R}^{n}$ are  vectors given by monomials of second and third order, respectively. 
We stress that the components of these vectors can be expressed using an indicial notation
\[
\bfG_i=\sum_{j,k=1}^{n}G_{ijk} D_j D_k, \quad   \bfH_i=\sum_{j,k,l=1}^{n}H_{ijkl}D_j D_k D_l, \qquad i=1,\ldots, n.
\] 

Equation \eqref{eq:PPV_d2} represents our high-fidelity, FOM which depends on the input parameters
$\omega,\bsn{\beta}$. 
The FOM can be solved in different ways and in the present work we consider two alternatives:  
a time marching scheme, i.e. a nonlinear Newmark algorithm, 
and an HB solver as developed in \cite{actuators21}.
It should be recalled that
in resonating MEMS an important output of interest is the FRF
in which a selected quantity, like the midspan deflection of a beam or the rotation
of a micromirror, is plotted versus $\omega$ for different $\bsn{\beta}$.
Indeed, the focus is on frequency stability for the following main reason: 
resonators operate close to a reference frequency where the behavior 
should be predictable. For instance,  in micromirrors the stability of the motion is required to guarantee
the performance during the line scanning process 
and predicting correctly the hardening and softening behavior is of paramount importance. 
As a consequence we are interested  in the steady state response of the device.
This is the direct output of HB approaches which express the solution as the sum of Fourier series.
However, HB solvers are not standard in commercial codes and might not be easily accessible.
Moreover, their cost rapidly increases with the size of the Fourier basis thus requiring dedicated
computing facilities.
On the contrary, time marching schemes are always available, but transients 
before reaching the steady state condition are 
often prohibitively long due to the large quality factors of MEMS.
As a consequence the choice of the solver is in general a trade-off which
strongly depends on the application at hand. Several examples are commented 
in Section~\ref{sec:mech} where details on the simulation settings are provided.

\subsection{Reduced order modeling through POD}

The first step in the construction of a POD-Galerkin ROM requires to generate 
a  matrix $\bfX\in \mathbb{R}^{n\times m}$, whose $m$ columns collect snapshots 
 of the FOM solutions, obtained for different values of the parameters $\omega, \bsn{\beta}$. 
If the FOM is solved by means of an HB approach, the snapshots for a given frequency are taken at regular intervals 
over one single period of the steady state response by reconstructing the displacement field starting from the Fourier coefficients. Otherwise, if time marching schemes
are employed, several alternatives are indeed available according to whether snapshots are taken 
in a condition close to the steady state or not. The influence of these choices is extensively 
investigated   in Section~\ref{sec:mech}.


Next, the  Singular Value Decomposition (SVD) of the matrix $\bfX$ is computed,
\[
\bfX = \bfU \bfSigma \bfV^T 
\]
where the columns of the orthonormal matrix $\bfU\in \mathbb{R}^{m\times m}$ are 
the left singular vectors, often called Proper Orthogonal Modes (POMs) in the literature 
\cite{kerschen2005,alfiobook,lu2019};
the columns of the orthonormal matrix $\bfV\in \mathbb{R}^{n\times n}$ are the right singular vectors. 
The diagonal elements of $\bfSigma\in \mathbb{R}^{m\times n}$ are the singular values of the matrix $\bfX$ and are conventionally ordered 
from the largest to the smallest. In particular, the rank of $\bfX$ is equal to the number of nonzero singular values, and  the optimal rank-$p$ approximation $\tilde{\bfX}$ of $\bfX$, in a least squares
sense, is given by the rank-$p$ SVD truncation 
\[
\tilde{\bfX}=\sum_{i=1}^p \sigma_i \bfU_i \bfV^T_i
\]
in which $\sigma_i$ is the  $i$-th singular value contained in the diagonal of  $\bfSigma$, 
$\bfU_i$ is the $i$-th column of $\bfU$ and $\bfV_i$ is the $i$-th column of $\bfV$. See, e.g., \cite{alfiobook} for further details. The SVD of $\bfX$ provides important insight into the energy distribution of the snapshots where the energy  of $\bfX$ is defined via the Frobenius norm:
\[
\veps(\bfX) = \|\bfX\|_F^2  =\sum_{i=1}^{n}\sum_{j=1}^{m} x^2_{ij} = \sum_{k=1}^{\min(m,n)} \sigma_k^2. 
\]
The POD approach selects the first $p$ most energetic POMs to build the POD-Galerkin ROM approximation:
\begin{equation}\label{eq:trialspace}
\bfD \approx 
\sum_{i=1}^p Q_i \bfU_i 
\end{equation}
where $Q_i$ are the ROM generalized coordinates. Hence, the POD-Galerkin approximation of 
$\bfD$ in \eqref{eq:trialspace} is given by a linear combination of POD modes, and the resulting trial subspace is optimal in the sense that it captures the highest possible energy content among all possible linear subspaces for any prescribed dimension $p$. Furthermore, the error in the snapshots approximation is related to the sum of the square of the singular values associated to the nonretained modes \cite{alfiobook}.

Once the linear trial POD subspace has been obtained, projecting the FOM \eqref{eq:PPV_d2} onto the POD subspace yields the structural dynamics geometric POD-Galerkin ROM, under the form of a $p$-dimensional nonlinear ODE system, whose solution provides the dynamics of the generalized coordinates: 
\begin{equation}\label{eq:PPV_POD}
\begin{split}
\bfM\aPOD \ddot{\bfQ} + \bfC\aPOD\dot{\bfQ} & + \bfK\aPOD\bfQ + \bfG\aPOD(\bfQ,\bfQ) \\
& +\bfH\aPOD(\bfQ,\bfQ,\bfQ) = \bfF\aPOD(\bfQ,\bsn{\beta},\omega,t), \quad t \in (0,T)
\end{split}
\end{equation}
where
\[
\bfM\aPOD=\bfU^T\bfM\bfU, \quad \bfC\aPOD=\bfU^T\bfC\bfU, \quad \bfK\aPOD=\bfU^T\bfK\bfU,
\]
\[
\bfF\aPOD=\bfU^T\bfF, \quad  G\aPOD_i=g\aPOD_{ijk} Q_j Q_k, \quad  H\aPOD_i=h\aPOD_{ijkl} Q_j Q_k Q_l,
\] 
with $\bfM\aPOD,\bfC\aPOD,\bfK\aPOD \in \mathbb{R}^{p\!\times\!p}$.
The computation of the vectors $\bfG\aPOD$ and $\bfH\aPOD$ entails $O(p^3)$ and $O(p^4)$ terms, respectively.
Note that the coefficients $g\aPOD_{ijk}$ and $h\aPOD_{ijkl}$ can be precomputed, and that the reduced problem can be assembled efficiently thanks to its polynomial nature, 
thus avoiding the use of  hyper-reduction techniques such as the (discrete) 
empirical interpolation method \cite{barrault2004anempirical,chaturantabut2010nonlinear,maday2008ageneral}.

\subsection{Solution of the Reduced Order Model}
\label{sec:solROM}

One of the greatest benefits of generating a POD-Galerkin ROM as the one in eq.\eqref{eq:PPV_POD} 
is the possibility to compute directly periodic solutions and trace the full FRF,
with both stable and unstable branches, by resorting to continuation codes 
either based on HB techniques or collocation approaches. 
Some well-known packages, suitable for small scale problems, are available in the literature. One of the most relevant examples is \texttt{Auto07p} \cite{Doedel}, a package that uses collocation methods in 	\texttt{FORTRAN} to perform numerical continuation and bifurcation analysis. 
Among other tools we can mention \texttt{Manlab}, a \texttt{Matlab} package that uses HB method and Asymptotic Numerical Method \cite{Guillot1,Guillot2}; \texttt{Nvlib} that also exploits HB methods \cite{NVLIB}; 
\texttt{COCO} that implements collocation methods and algorithms for bifurcation detection \cite{COCO}. 
Another excellent package able to perform the continuation of ODEs is \texttt{BifurcationKit} \cite{BifurcationKit}, an emerging toolkit for Julia language that provides continuation methods for ODEs and PDEs. These packages usually provide the ability to distinguish between stable and unstable branches, locate bifurcation points 
and follow alternative branches of the solution. 
We highlight that the same versatility is difficult to achieve with a FOM. Indeed, even if a HB formulation with continuation has been recently proposed in \cite{KerschenHB2015,actuators21} for large scale problems, 
computing times are not compatible with their application at the design or prototyping levels.

In this work, we compute solutions with \texttt{Manlab}, which has an impressive capability 
to perform accurate bifurcation analysis and to exploit the Asymptotic Expansion Method \cite{cochelin2009high} by ensuring a good balance between computational time and accuracy, 
provided that the problem can be re-written in quadratic form.

\subsection{POMs and reconstruction vectors in NNM}
\label{sec:NNM}

An important feature of an effective ROM is the capability to identify an invariant subspace for the system dynamics, i.e. trajectories initiated along the subspace remain within the subspace itself in the full order solution. 
In linear systems each mode defines an invariant plane in the phase space, hence linear projection methods as the modal decomposition provide an excellent tool for generating ROMs.
On the other hand, in presence of geometric nonlinearities, invariance of modal subspaces is not guaranteed as underlined in past works, for instance by Amabili and Touz\'e \cite{amabili07} and Haller \cite{buza2021using}. Indeed, the invariant manifold tangent at the origin to a given modal subspace is a curved hypersurface that requires nonlinear projection methods. 
In this framework, the parametrization method initially formulated by Haro and De la Lave \cite{haro2006parameterization1,haro2006parameterization2,haro2016parameterization} was recently applied to large scale finite element systems of mechanical structures \cite{vizza2020,nld21,jain2021compute} thanks to the Direct Parametrization approach.  
The fundamental idea of this class of methods is to parametrize the dynamics of the system along the invariant manifold associated to one of its eigenfunctions. This requires the introduction of a nonlinear change of coordinates between nodal displacements and the parametrization coordinates. Using the formulation proposed in \cite{vizza2020}, the nonlinear coordinates change for a single master-mode reduction in an undamped mechanical system is expressed as:
\begin{align}
	\label{eq:reconstr}
	\bfD = &\bfphi_{m}R + \hat{\bfa}R^{2} + \hat{\bfb}S^{2} + \hat{\bfc} R^{3} + \hat{\bfu} RS^{2} + O(|R,S|^{4}),
\end{align}
where $\bfphi_{m}$ denotes the  eigenmode associated to the master mode. Here $\hat{\bfa}$, $\hat{\bfb}$, $\hat{\bfc}$, and $\hat{\bfu}$ are higher order reconstruction vectors used to map the parametrization coordinates $R$ and $S$ to the physical displacement $\bfD$. 
The reconstruction vectors in eq.\eqref{eq:reconstr} apply a correction with respect to a simple modal decomposition approach by accounting for the coupling between master and slave modes. As shown by Buza \cite{buza2021using}, the projection of the eigenfunctions of the system along quadratic reconstruction vectors provides a solid framework to identify the modes that better describe the curvature of the invariant manifold. 
This last result is the natural extension of what remarked by Amabili and Touz\'e \cite{amabili07}
where the trial space identified by POD was interpreted as the best linear approximation of the nonlinear normal mode. The consequence is that, in order to properly build a ROM relying on  methods such as POD, one needs to introduce also bases that allow a correct approximation of the manifold curvature. This is highlighted in the results section, where qualitative changes in the predicted structural response are obtained by adding PODs with apparent negligible energy contribution, 	however showing a high curvature related to the invariant manifold of the system. This effect can also be observed in their shape, which resembles that of the reconstruction vectors provided by the DP.

\section{Purely mechanical applications}
\label{sec:mech}

Four mechanical benchmark cases are here proposed
to discuss the accuracy of POD-Galerkin ROMs: a doubly clamped beam, two micromirrors and a shallow arch.
In all cases, FEM meshes are made of wedge quadratic elements (``extruded'' isoparametric elements with 15 nodes). 
For each example we report in \ref{sec:Annex} the computation time required by the FOM, the offline and the online stages of the ROM.

\subsection{Doubly clamped beam}
\label{sec:ccbeam}

\begin{figure}[ht]
\centering
\includegraphics[width = .85\linewidth]{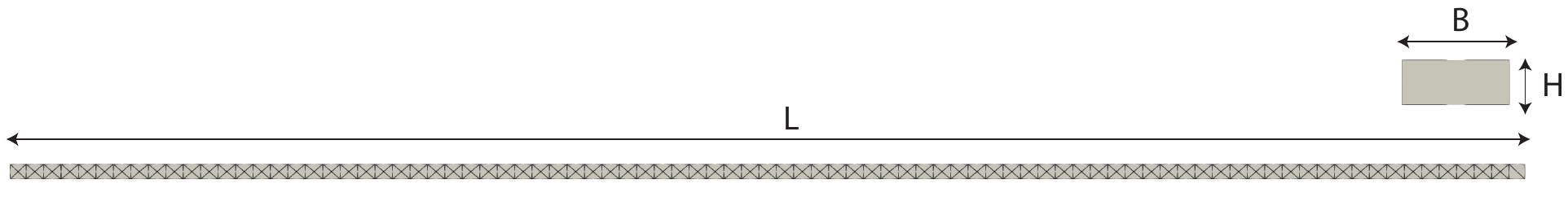}
\caption{Geometry and mesh of the doubly clamped beam, front view and cross-section}
\label{fig:beam}
\end{figure}

Let us consider a doubly clamped beam of length $L=1000$\micr with a rectangular cross-section
of dimensions 10\micr$\times$24\micr, as depicted in Fig.\ref{fig:beam}. 

This academic example simulates realistic MEMS resonators like those analysed
in \cite{jmems20reso}.
A rather coarse mesh with 2607 nodes has been employed. 
Indeed, this example is used to discuss extensively the features of the ROM and 
every result is compared with 
reference FOM solutions that are affordable only with relatively coarse meshes.
The beam is made of isotropic polysilicon
\cite{jmems04}, with density $\rho=2330$\,Kg/m$^3$, Young modulus $E=167$\,GPa and 
Poisson coefficient $\nu=0.22$. We select a fixed quality factor $Q=50$. 
The device vibrates according to its first bending mode at 
$f_0= 87141$\,Hz. The first five eigenfrequencies are reported in Table~\ref{tab:CC_freq}.

\begin{table}[h!]
\centering
\begin{tabular}{ |c|c|c|c|c|c| } 
 \hline
Eigenmode & 1 & 2 & 3 & 4 & 5\\ 
 \hline
Frequency  [kHz]& 87.141 & 208.45 & 240.03 & 470.10 & 572.18 \\ 
 \hline
\end{tabular}
 \caption{Doubly clamped beam: eigenfrequencies}
\label{tab:CC_freq}
\end{table}

\begin{figure}[h!]
\centering
\includegraphics[width = .85\linewidth]{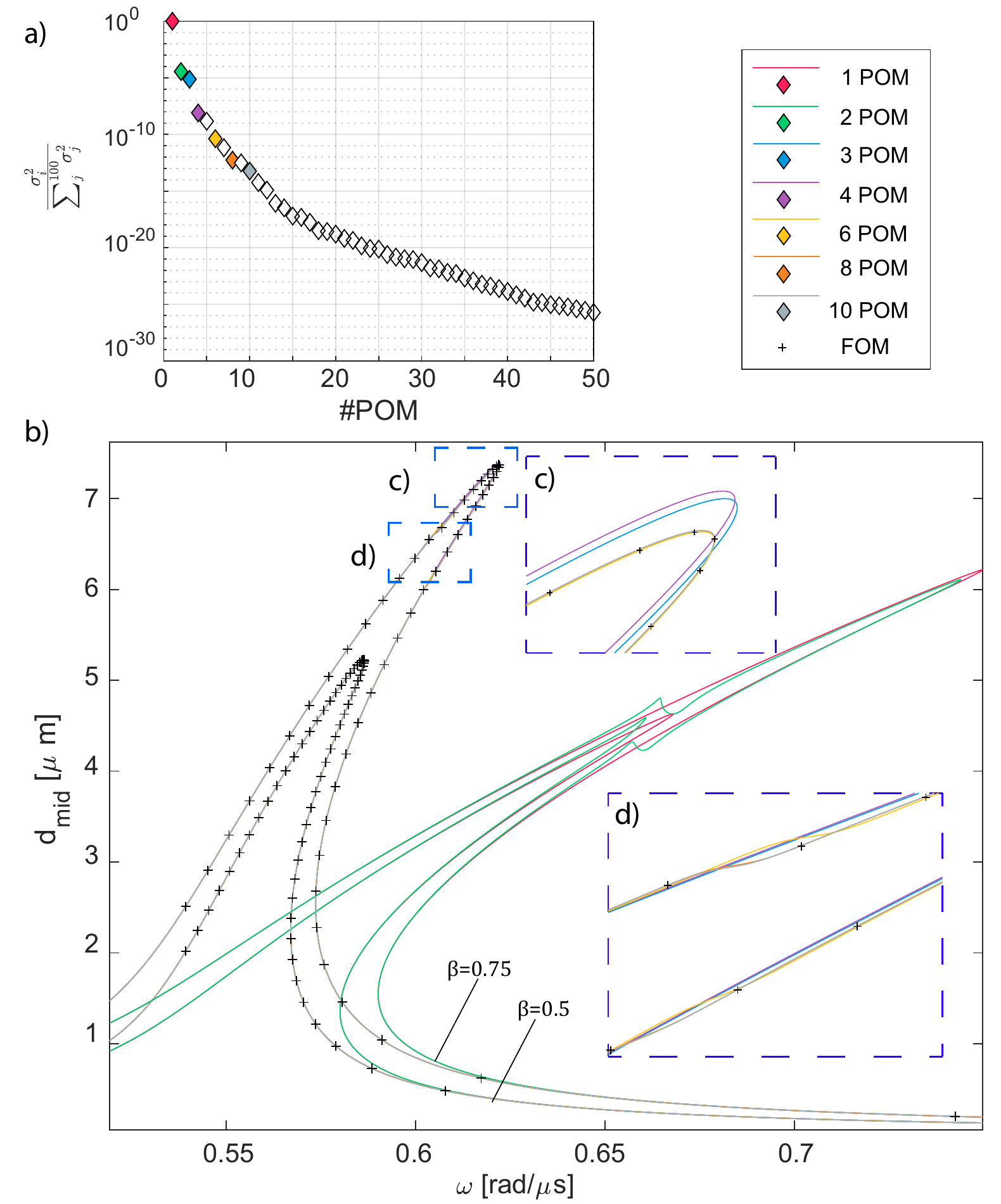}
\caption{Doubly clamped beam: convergence with respect to the number of POMs retained in the ROM. Figure a): relative energy content of POMs. Figure b): FRF computed with each ROM subspace considered. Figures c) and d): enlarged views of the resonance peak and the internal resonance interaction region, respectively.}
\label{fig:FRF_NPOM_HB_beam}
\end{figure}
The external excitation is provided by a body load proportional to the first eigenmode 
$\bfF=\bfM \bfphi_1 \beta \cos(\omega t)$ with $\beta$ load multiplier.

\begin{figure}[h!]
\centering
\includegraphics[width = .8\linewidth]{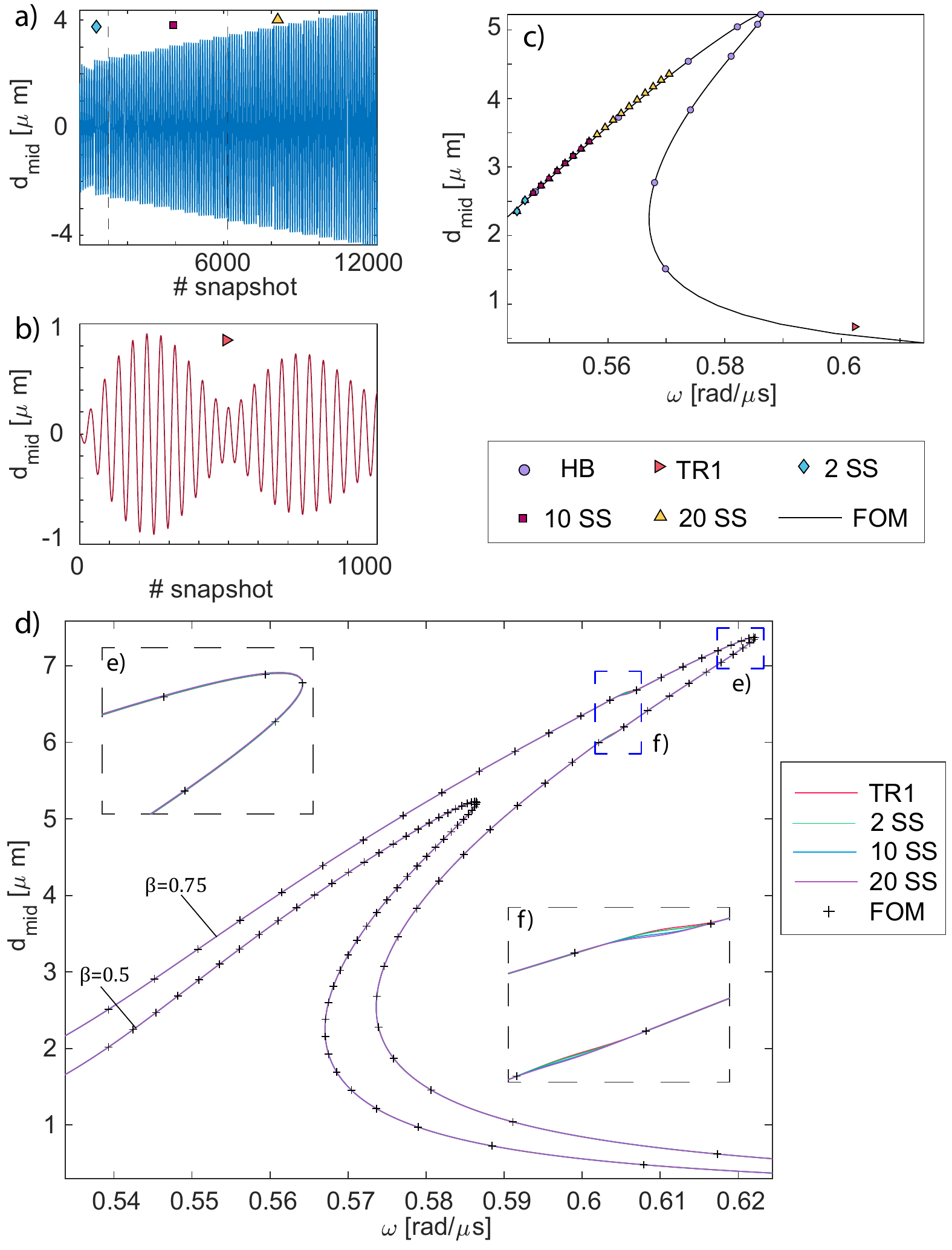}
\caption{Doubly clamped beam. Comparison between training with time-marching or HB snapshots. 
Figure a): time histories of the midspan displacement corresponding to the snapshots set collected close to Steady State (SS). The frequency is swept upwards and the data are collected after a fixed number of time steps. The time histories are simulated sequentially and jumps denote a
change of the forcing frequency. Figure b): time history of the midspan displacement corresponding to the snapshot set collected in a transient (TR) case. Figure c): sampling frequencies of the datasets  HB, SS and TR. The SS datasets differ due to the number of frequencies sampled. 
Figure d) presents the FRFs computed with the ROMs. Figure e) and f): close-up of the resonance peak and of the internal resonance region.}
\label{fig:beam_samples}
\end{figure}

The training data can be generated with different methods and a variable 
number of snapshots. Since we aim at modeling the steady-state response of the system, 
HB solutions are the ideal candidates to generate representative data. 
In a first application, we consider a training dataset 
computed with HB and $\beta=0.5\,\mu$N
and consisting of a total of 500 snapshots generated from one period of 10 frequency samples along the FRF, 
represented by the violet circle markers in Figure \ref{fig:beam_samples}c.

First, we address the convergence with respect to the number of POMs retained in the subspace.
The POD method usually adopts the relative energy as a convergence measure to select the space dimension.
The relative energy content of each POM is depicted in Figure~\ref{fig:FRF_NPOM_HB_beam}a. 
It is worth stressing that the first POM represents 99.995 $\%$ of the energy and the second POM contains only 0.0035 $\%$ of the energy.
To test the convergence of the ROM, we consider seven different subspaces with 1,2,3,4,6,8 and 10 POMs, respectively. The resulting ROMs are tested on 
$\beta=0.5\,\mu$N (i.e.\ the same forcing level seen during the training stage) and $\beta=0.75\,\mu$N. 
The resulting FRFs are plotted in Figure~\ref{fig:FRF_NPOM_HB_beam}b.
We notice that trial spaces with less than 3 POMs are inadequate to describe the dynamics,
while richer spaces provide a very good accuracy. 
It should be remarked also that increasing the number of POMs a mild 
1:5 internal resonance is evidenced, as predicted by the high-fidelity FOM
(see Figure \ref{fig:beam_samples}d).
 
\begin{figure}[h!]
\centering
\includegraphics[width = .65\linewidth]{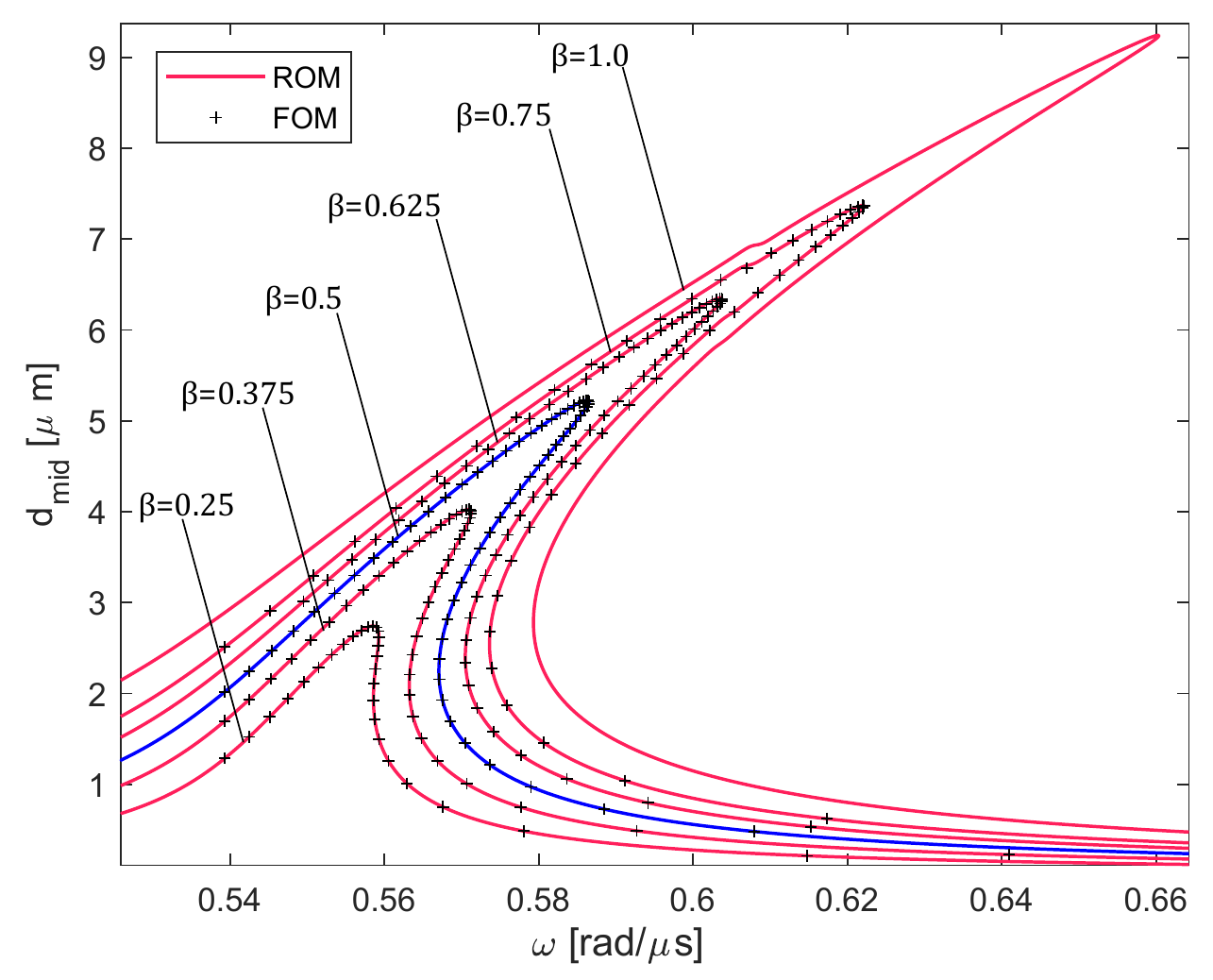}
\caption{Doubly clamped beam: FRFs computed with the ROMs with 6 POMs compared with the FOM solution. 
The blue line denotes the training curve, while the red ones are 
solutions computed for different levels of the forcing. The cross markers represent the FOM solutions}
\label{fig:FRF_beam}
\end{figure}

Next, keeping 6 POMs in the trial space, a choice that guarantees a good balance 
between efficiency and accuracy,
the ROM is tested on different forcing levels 
$\beta=0.25,0.375,0.45,0.5,0.625,0.75\,\mu$N. 
The results are reported in Figure~\ref{fig:FRF_beam}, always compared to  the solutions
of the HB-FOM.
The two families of simulations are almost exactly superimposed, hence
proving the high predictive capability of POD,
well beyond the training range. Also the correct evolution of the unstable branch
is reproduced, from mildly to strongly hardening at increasing actuation levels.

Nevertheless, HB solutions might not be accessible in general (e.g.\ in commercial codes), 
or might be too costly to generate. In these cases time-marching methods are the only option available to generate snapshots.
To highlight  the possible differences with HB-FOM solutions 
in this simple and small example, 
we consider four datasets computed with time marching methods and $\beta=0.5\,\mu$N: 
1) 1242 snapshots generated from a response close to the steady state (SS) at 2 different frequencies; 
2) 6210 snapshots generated from a response close to SS at 10 different frequencies  (including the ones of set 1); 
3) 12420 snapshots generated from a response close to SS at 20 different frequencies  (including the ones of set 2); 
4) 1000 snapshots generated from a fully transient (TR) response at one single frequency.

The sampling points and  plots of the time marching datasets are reported in 
Figures~\ref{fig:beam_samples}a, \ref{fig:beam_samples}b and \ref{fig:beam_samples}c.
The corresponding ROMs with 6 POMs are displayed in Figure~\ref{fig:beam_samples}d. 
We notice that all the solutions are almost exactly superposed and 
the localized 1:5 internal resonance is the only portion of the FRFs where minor differences can be appreciated (Figure~\ref{fig:beam_samples}f).
In particular, the capability of the technique to predict steady state solutions starting from
fully transient data (case 4) is impressive and very promising for MEMS applications
where large quality factors generally prevent time marching schemes from reaching steady state conditions
within reasonable time frames.

\begin{figure}[h!]
\centering
\includegraphics[width = .7\linewidth]{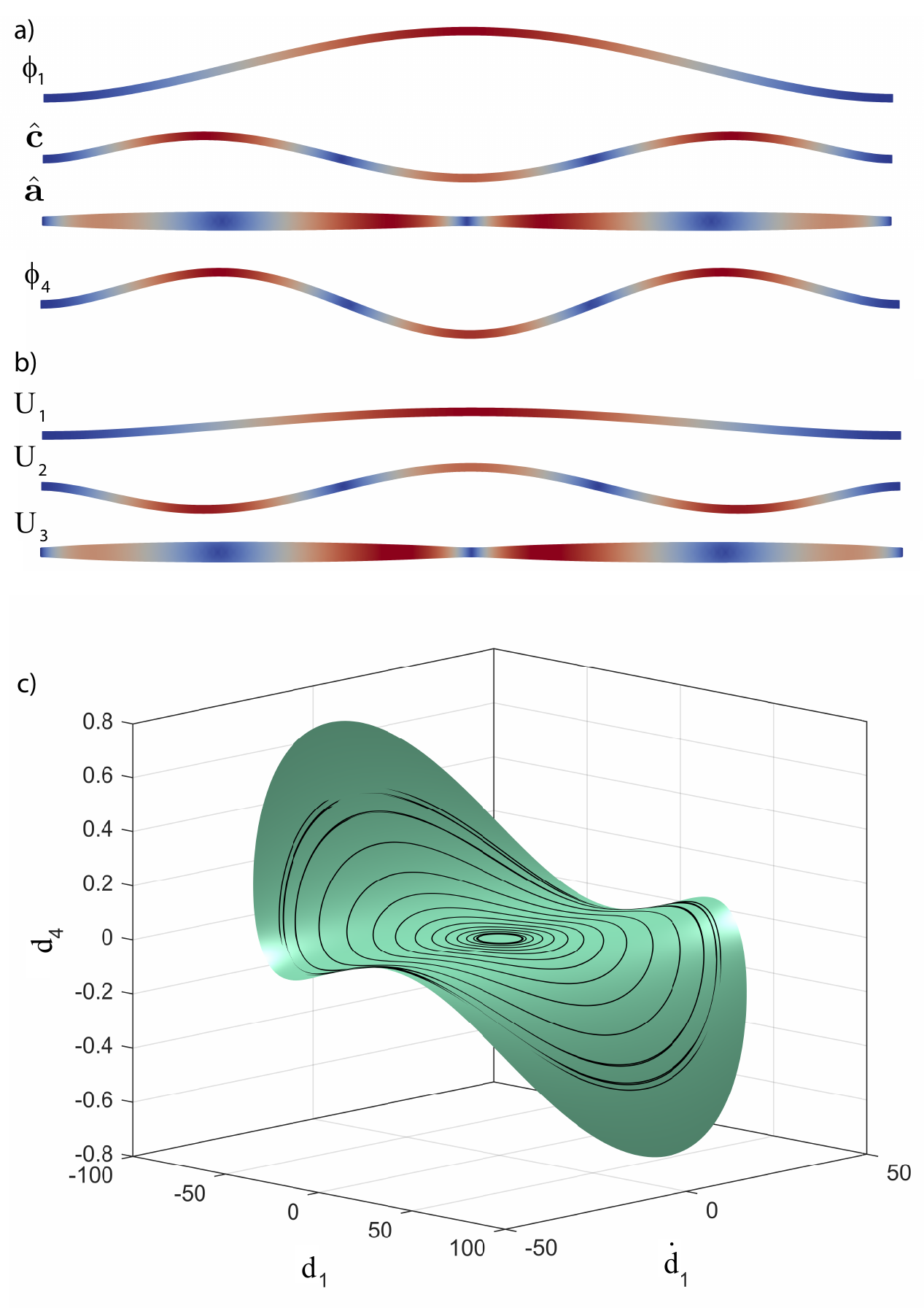}
\caption{Doubly clamped beam. Comparison between the POD subspace bases and DP. Figure a) shows the parametrization up to cubic order and the third eigenmode. Figure b): first three POMs. Figure c): invariant manifold computed with DP (green surface) and  orbits computed with POD (black lines). The manifold is defined in the phase space composed by the first eigenmode displacement and velocity and the fourth eigenmode displacement.}
\label{fig:beam_manifold}
\end{figure}

\subsubsection{Connection with modal methods and DP}

As put in evidence in Figure~\ref{fig:FRF_beam},
a major improvement comes from the inclusion of the third POM. 
The first three POMs are represented in Figure~\ref{fig:beam_manifold}b. 
The difference between the trial spaces  with 2 and 3 POMs
is not easy to appreciate through a direct application of the 
energy criterion. 
From a physical perspective, the first two POMs closely resemble 
the first and the second symmetric bending modes, 
while the third POM is related to high frequency axial/contraction modes.

In linear methods based on modal subspaces, 
see e.g.\ \cite{givois,cyrilstep},
it is now acknowledged that this type of modes 
must be imperatively included in the selected subspace to guarantee convergence,
but they are difficult to identify a priori.
These modes indeed provide an important correction to the stress field that 
cures the over-hardening typical of linear techniques.
The automatic identification of such a contribution can be considered as a major benefit of the POD
over modal techniques.

Considering now the parametrization methods discussed in Section~\ref{sec:NNM},
the reconstruction vectors of eq.\eqref{eq:reconstr} 
are plotted in Figure~\ref{fig:beam_manifold}a.
We start noticing the striking analogy between the first eigenmode $\bfphi_1$, 
the quadratic displacement-dependent term $\hat{\bfa}$, 
the cubic displacement-dependent term $\hat{\bfc}$ and the first three POMs.
 
Next, in Figure~\ref{fig:beam_manifold} the manifold of the DP is 
compared with orbits computed through POD by considering the phase space 
composed by the projection of the displacement and velocity on 
the first eigenmode and the projection of the displacement on the fourth eigenmode.
The ROM solution lies almost perfectly on the approximated manifold, thus showing 
that the subspace computed with data-driven methods converges to the one computed 
with asymptotic expansions.

\subsection{Micromirrors}
\label{sec:mirrors}

Scanning micromirrors are witnessing explosive growth in recent years
due to successful applications ranging from pico
projectors for Augmented Reality (AR) lenses \cite{ARlenses}, to 3D scanners
for Light Detection and Ranging (LiDAR) application.
In this section we address two micromirrors with different nonlinear behavior.
These devices are intrinsically nonlinear due to inertia effects associated with  large rotations, and they present a frequency-amplitude dependence which might be 
either hardening or softening according to the specific layout. 
The correct quantitative prediction of nonlinear effects 
is a tough benchmark for any FOM or ROM.
Recently, the authors have developed a large scale HB approach 
in \cite{actuators21} for the analysis of piezo mirrors which is here utilized as FOM
to generate snapshots.
It is worth stressing that classical ROM techniques like the Implicit Condensation \cite{ijnm19} and
Modal Derivatives \cite{md2016} fail to provide the required accuracy. 

\subsubsection{Micromirror 1}
\label{sec:mirror1}

\begin{figure}[h]
\centering
\includegraphics[width = .7\linewidth]{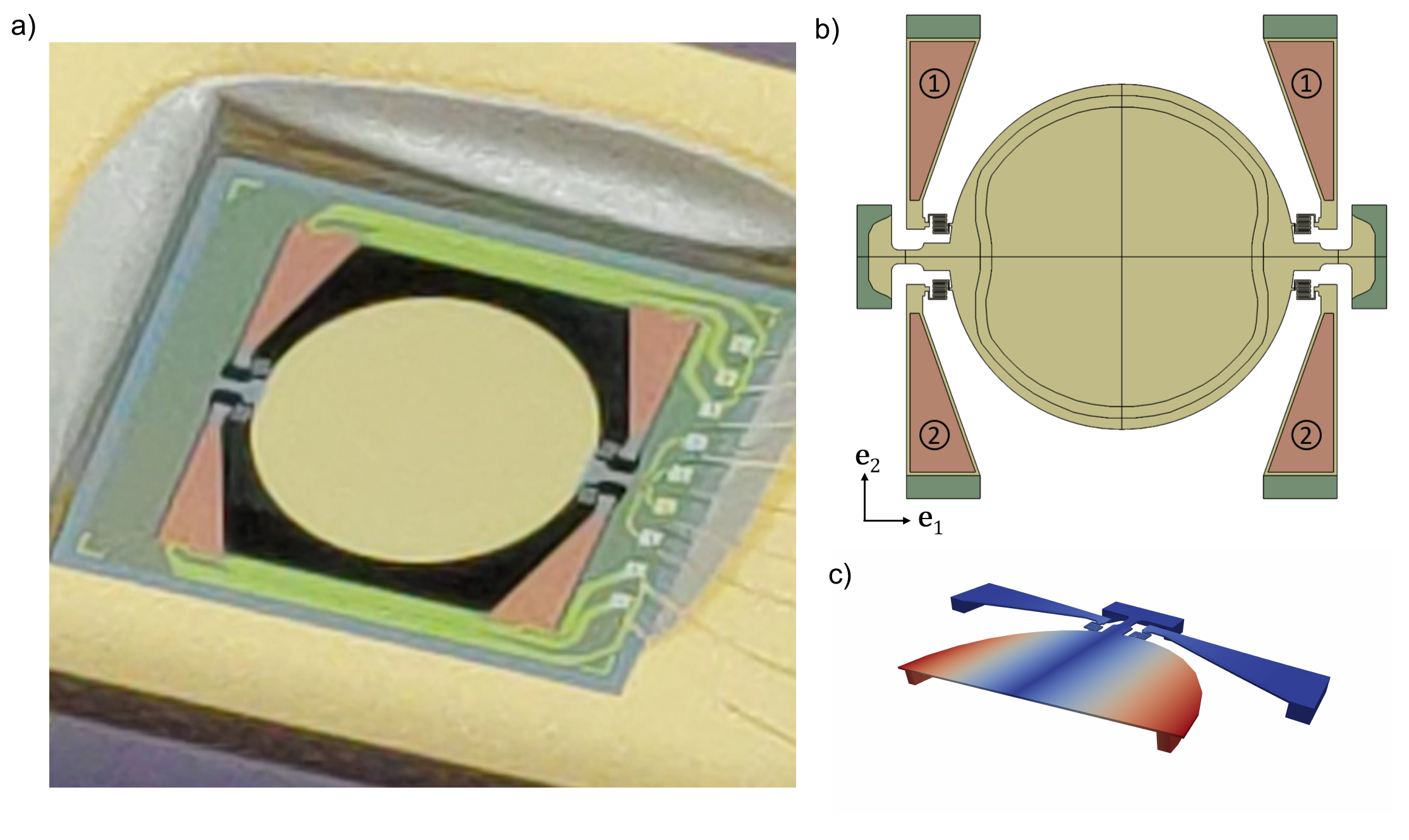}
\caption{Micromirror 1. Figure a): photo of the real device. Figure b): top view of the layout. Piezoelectric patches are in light brown and the numbers characterize the actuation scheme \cite{actuators21,jmemsmirror}. Figure c): first torsional 
mode. The eigenmode consists of a rotation of the micromirror plate as shown by the color-map of the displacement magnitude.}
\label{fig:mirror}
\end{figure}

An optical image of the first micromirror, fabricated by STMicroelectronics
\cite{actuators21,jmemsmirror}, 
is presented in Figure~\ref{fig:mirror}a. The top view is reported in Figure~\ref{fig:mirror}b. The reference FOM is built considering only half of the structure to exploit symmetry as illustrated in Figure~\ref{fig:mirror}c. 
The central circular reflecting surface is directly attached to the substrate
with two short torsional beams (springs), while the rotation of the mirror is induced
by trapezoidal beams which are connected to the mirror
through folded compliant springs. The beams are actuated by piezo-patches, i.e.\ PZT layers with a thickness of 2\micr appearing in Figure~\ref{fig:mirror}b as light brown areas.
The mirror disk has a diameter of 3000\micr and the lower
surface has been reinforced with a curvilinear
support in order to minimize the dynamic deformation
of the mirror itself.

The mirror is made of single crystal silicon with [100] orientation \cite{hopcroft2010} and 
the first five eigenfrequencies are reported in Table~\ref{tab:mirr1_freq}.

\begin{table}[h]
\centering
\begin{tabular}{ |c|c|c|c|c|c| } 
 \hline
Eigenmode & 1 & 2 & 3 & 4 & 5\\ 
 \hline
Frequency  [kHz]& 2.258 & 7.238 & 23.378 & 23.426 & 56.046 \\ 
 \hline
\end{tabular}
 \caption{Micromirror 1: eigenfrequencies}
\label{tab:mirr1_freq}
\end{table}

\begin{figure}[h!]
\centering
\includegraphics[width = .8\linewidth]{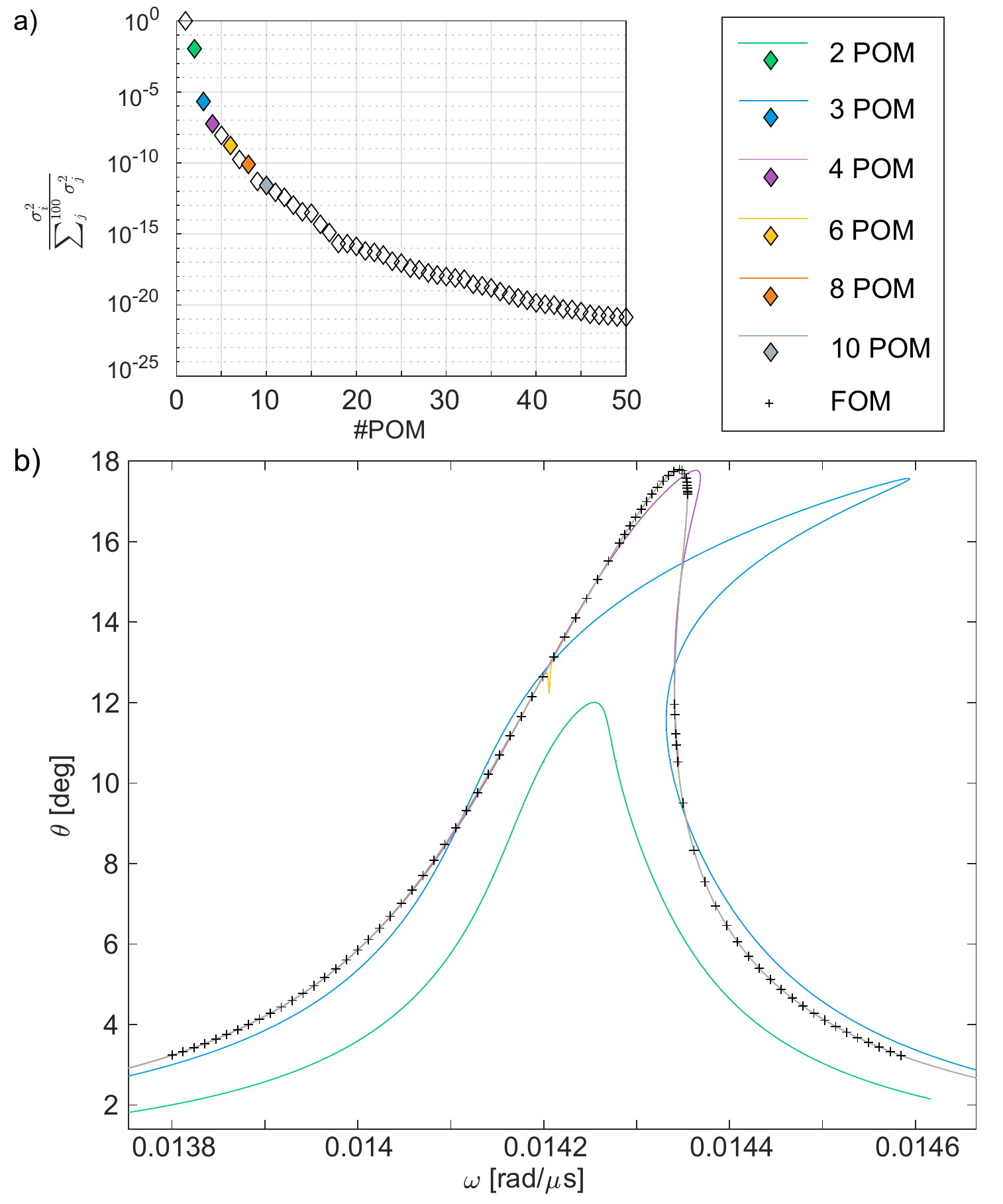}
\caption{Micromirror 1: convergence of the ROM. Figure a): energy of each POM. 
Figure b): FRF for the ROMs including an increasing number of POM. 
The POD is tested on the training set. 
The ROM simulations progressively converge to the correct FOM solution. 
With 6 POMs a small resonance effect occurs close to $\omega=0.0142 [\text{rad}/\mu s]$ 
that is eliminated by further increasing the number of POMs}
\label{fig:lynx_enor}
\end{figure}

The quality factor is set to $Q=100$.
In our investigation, for the sake of simplicity, 
we replace the piezoelectric actuation with a body force proportional 
to the first eigenmode $\bfF(t)=\bfM  \bfphi_1  \beta\cos(\omega t)$ with $\beta$ load multiplier.
The FEM model in this benchmark contains a total of 15341 nodes.

In the training stage, in which the FOM has been solved with an HB approach, 5000 snapshots have been generated from frequency samples uniformly distributed over the FRF setting $\beta=0.3\,\mu$N.
Even if a POD-Galerkin ROM with 1 or 2 POMs collects more than $99.99\%$ of the energy content
(see Figure~\ref{fig:lynx_enor}a), 
the trial space should contain at least the first 6 POMs to achieve convergence. Nevertheless the ROM with 6 POMs displays a small resonance effect close to $\omega=0.0142$ that is eliminated by using 8 POMs.
This can be appreciated from Figure~\ref{fig:lynx_enor}b where several FRFs have been computed 
with an increasing number of POMs.

\begin{figure}[ht]
\centering
\includegraphics[width = .7\linewidth]{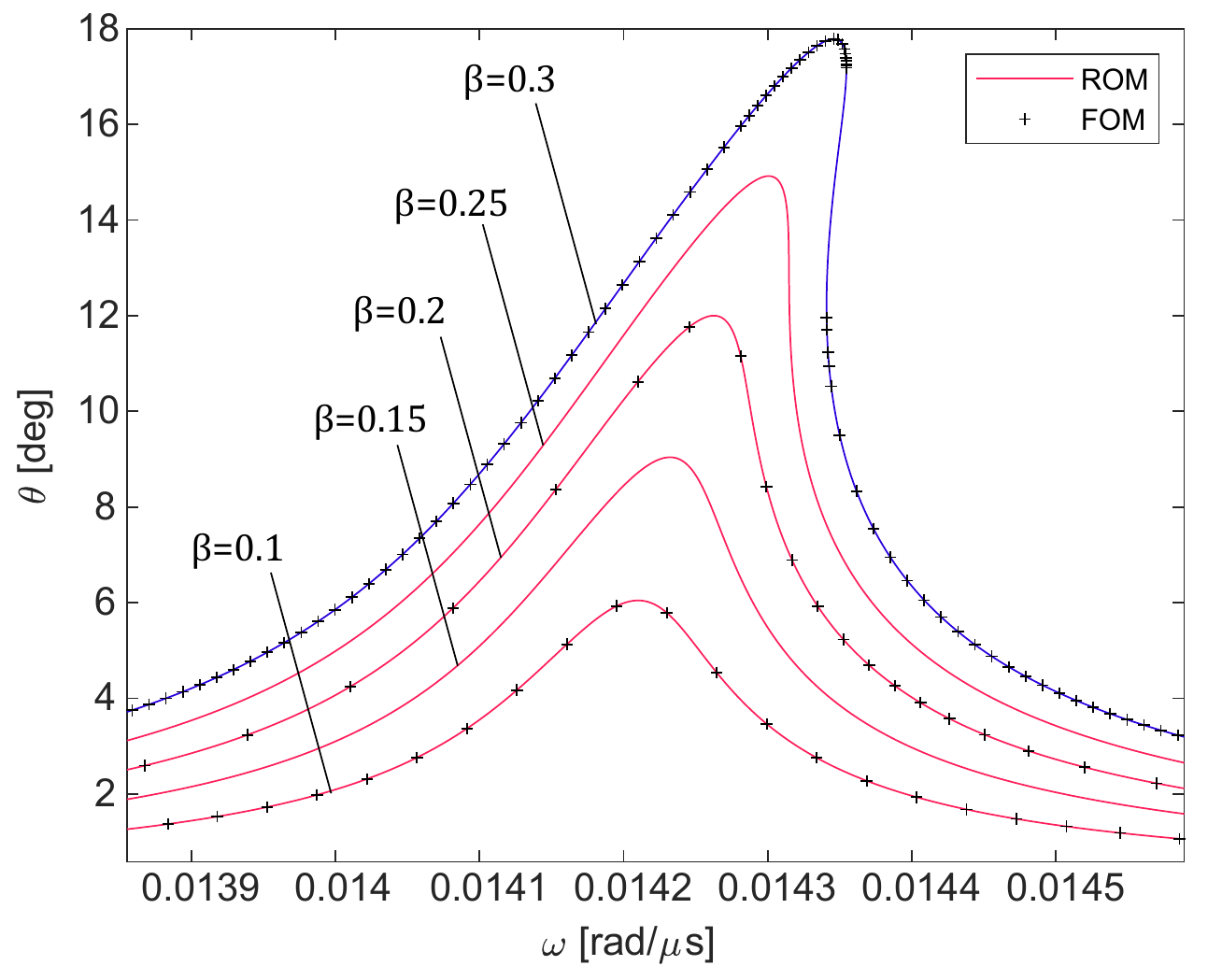}
\caption{Micromirror 1: FRFs computed with the ROMs with 10 POMs compared with the FOM solution. The red continuous lines represent the ROM solutions computed in conditions different from the training data, the blue one marks the training condition curve. The cross markers represent the FOM solutions}
\label{fig:lynxFRF}
\end{figure}

For the subsequent analyses, a trial space with 8 POMs is retained.
Different levels of the forcing have been tested as illustrated in Figure~\ref{fig:lynxFRF}
showing that an excellent agreement is achieved together with a remarkable
predictive capability. 
Only few curves are validated against the corresponding FOM solution because of the prohibitive computational cost entailed by the HB-FOM model. \\

\begin{figure}[h!]
\centering
\includegraphics[width = .7\linewidth]{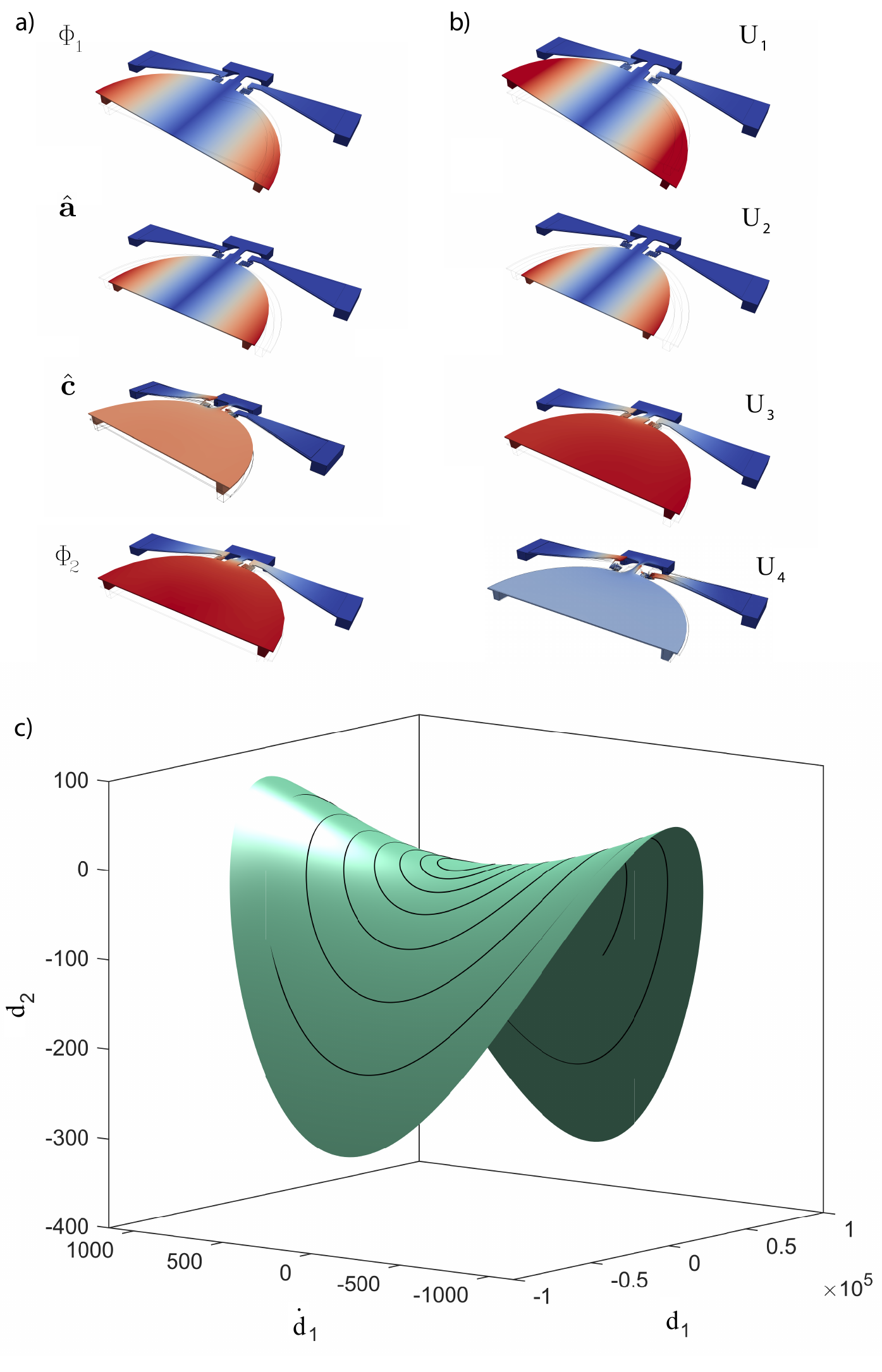}
\caption{Micromirror 1: comparison between the POD subspace bases and DP. Figure a): parametrization up to cubic order and the second eigenmode. Figure b): first three POMs. Figure 3): invariant manifold computed with DP (green surface) and orbits computed with POD (black lines). The manifold is defined in the phase space composed by the first eigenmode displacement and velocity and the second eigenmode displacement}
\label{fig:lynx_manifold}
\end{figure}

Similarly to the doubly clamped beam, some remarks concerning the POD subspace are worth stressing.
The first 3 POMs are depicted in Figure~\ref{fig:lynx_manifold}.
From a physical point of view, the first POM corresponds to the linear torsional eigenmode, 
while the second one corresponds to a contraction of the mirror plate 
that applies a correction to the linearized torsion. 
Indeed a linearized rotation, when extended to large angles, induces a non physical stretch 
of the structure.
The higher order POMs correspond to the membrane and axial deformation of the deformable springs
and beams.

Considering the parametrization of the DP, 
in Figure~\ref{fig:lynx_manifold}a we also plot the reconstruction vectors of
eq.\eqref{eq:reconstr}.
Apart from the obvious correspondence between the linearized mode and the first POM,
we notice a strong correspondence
between $\hat{\bfa}$ and the second POM  and between $\hat{\bfc}$ and a combination of the third and the fourth POMs.

Finally, Figure~\ref{fig:lynx_manifold}c presents the solution manifold 
in a phase space composed by the projection of the displacement and velocity on 
the first eigenmode and of the displacement on the second eigenmode.
The manifold obtained with the DP (continuous green surface) and the orbits  obtained through the POD-Galerkin ROM
are nearly superimposed as remarked also in the previous example.
The strong connections between the POD and the DP emerge hence as 
a distinctive feature of the POD approach.

\subsubsection{Micromirror 2}

\begin{figure}[h!]
\centering
\includegraphics[width = .85\linewidth]{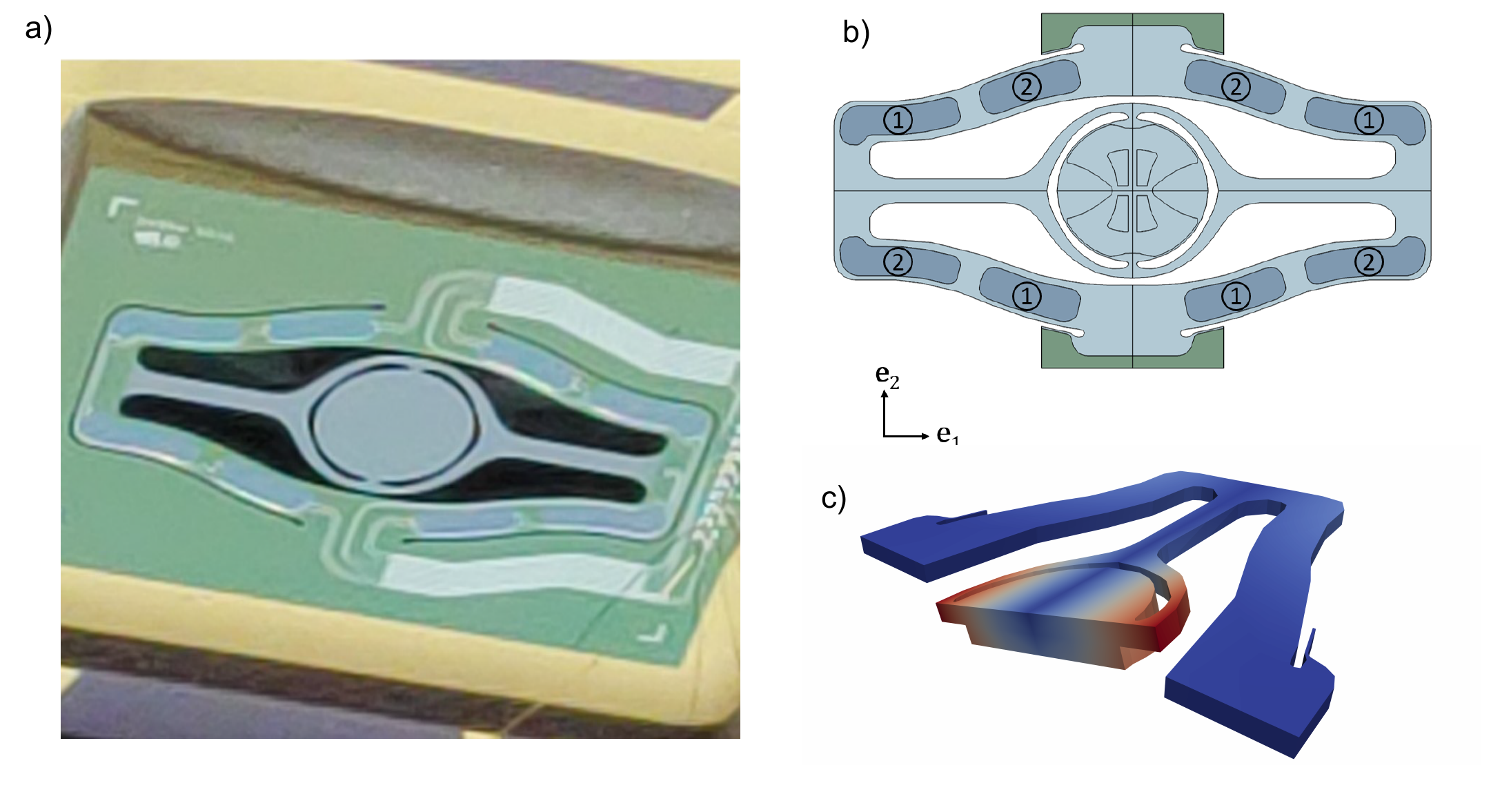}
\caption{Micromirror 2. Figure a): optical picture of the micromirror. 
Figure b): top view of the layout. Piezoelectric patches are in light brown and the numbers characterize the actuation scheme \cite{actuators21,jmemsmirror}. Figure c): 
third (torsional) eigenmode that is actuated during operations.}
\label{fig:perseus}
\end{figure}

The second mirror addressed, also fabricated by ST Microelectronics, is illustrated in Figure~\ref{fig:perseus}. In this case the mirror is suspended to a gimbal rather 
than being directly attached to the substrate with torsional springs. 
As a consequence, its nonlinear behavior turns softening. 

The frequency of the torsional mode is $29271$\,Hz 
and the quality factor has been set to $Q=1000$. The first 5 eigenfrequencies are listed in 
Table~\ref{tab:mirr2_freq}.
Also in this application, we replace the piezoelectric actuation method with a body force proportional to the third eigenmode $\bfF(t)=\bfM  \bfphi_3 \beta\cos(\omega t)$ with $\beta$ load multiplier.

\begin{table}[ht]
\centering
\begin{tabular}{ |c|c|c|c|c|c| } 
 \hline
Eigenmode & 1 & 2 & 3 & 4 & 5\\ 
 \hline
Frequency  [kHz]& 11.080 & 18.533 & 29.271 & 41.667 & 68.848 \\ 
 \hline
\end{tabular}
 \caption{Micromirror 2: eigenfrequencies}
\label{tab:mirr2_freq}
\end{table}

\begin{figure}[h!]
\centering
\includegraphics[width = 0.85\linewidth]{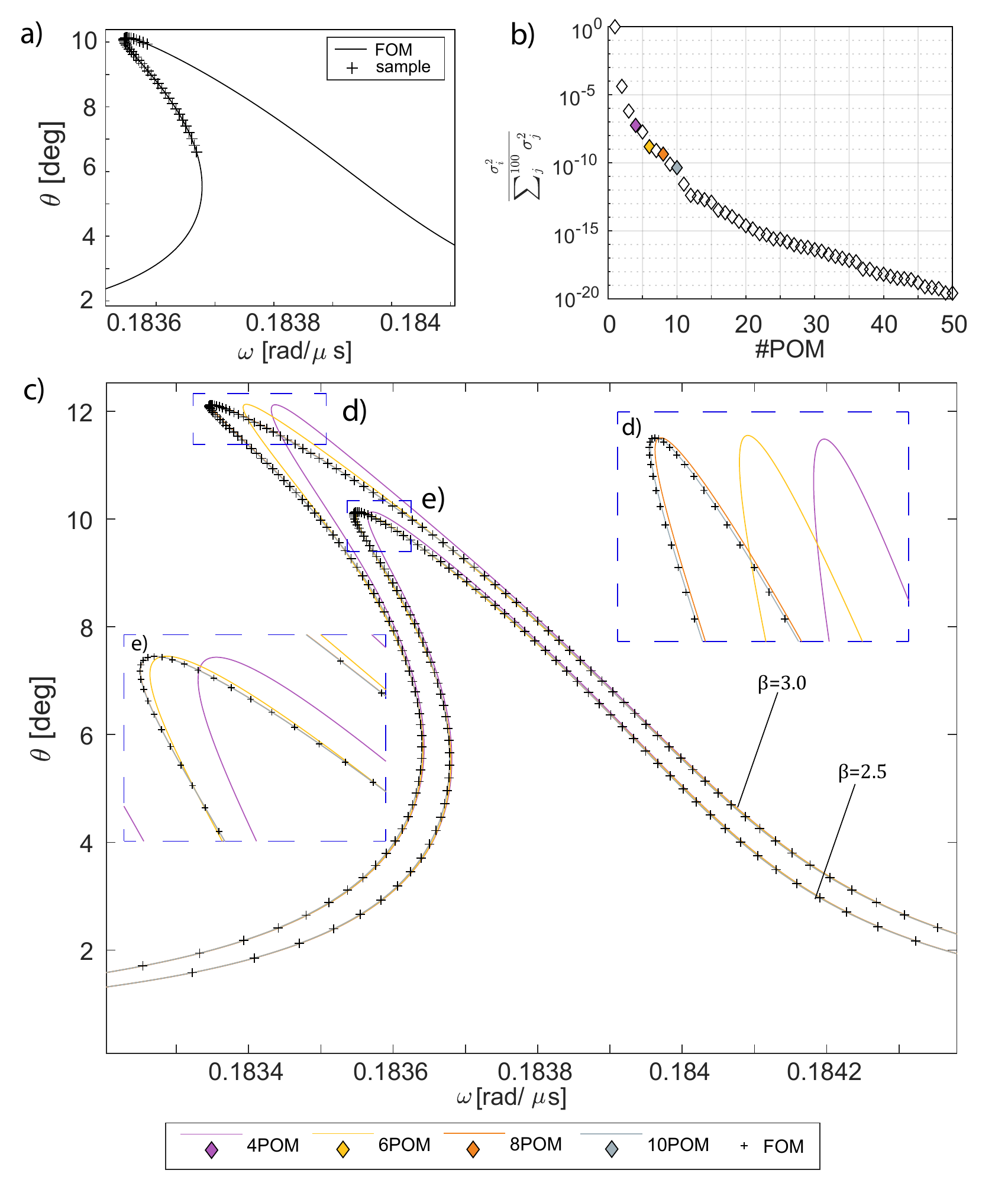}
\caption{Micromirror 2: convergence of the ROM.  Figure a): points of the FOM FRF 	utilized in the training phase to generate the snapshots. Figure b): energy content of each POM. Figure c): FRFs resulting from each ROM built with a increasing number of POMs. 
The POD is tested on the training set and on a second higher forcing level.}
\label{fig:FRF_POM_ok_perseus}
\end{figure}

This benchmark, which looks rather similar to the previous one, turns 
into a tough challenge for simulation approaches, mainly because the torsional mode is not the lowest-frequency one (it is the third) and is not well separated from the other modes.
Indeed, in this case even the DP technique requires a high order expansion and the
quadratic formulation in \cite{nld21} fails.
The POD, on the contrary, performs indistinctly well.
The training stage is performed considering $\beta=2.5\,\mu$N and generating a total of 2000 snapshots.
The distribution of the frequency samples on the FRF computed with the HB-FOM is presented in 
Figure~\ref{fig:FRF_POM_ok_perseus}a, while
Figure~\ref{fig:FRF_POM_ok_perseus}b collects the associated POMs.
Testing different ROMs built with an increasing number of bases we obtain the results displayed in Figure~\ref{fig:FRF_POM_ok_perseus}c. 
The convergence of the subspace is consistent with the increasing number of POMs, and this can be appreciated from the enlarged views in 
Figures~\ref{fig:FRF_POM_ok_perseus}d and \ref{fig:FRF_POM_ok_perseus}e. 

\begin{figure}[h!]
\centering
\includegraphics[width = .7\linewidth]{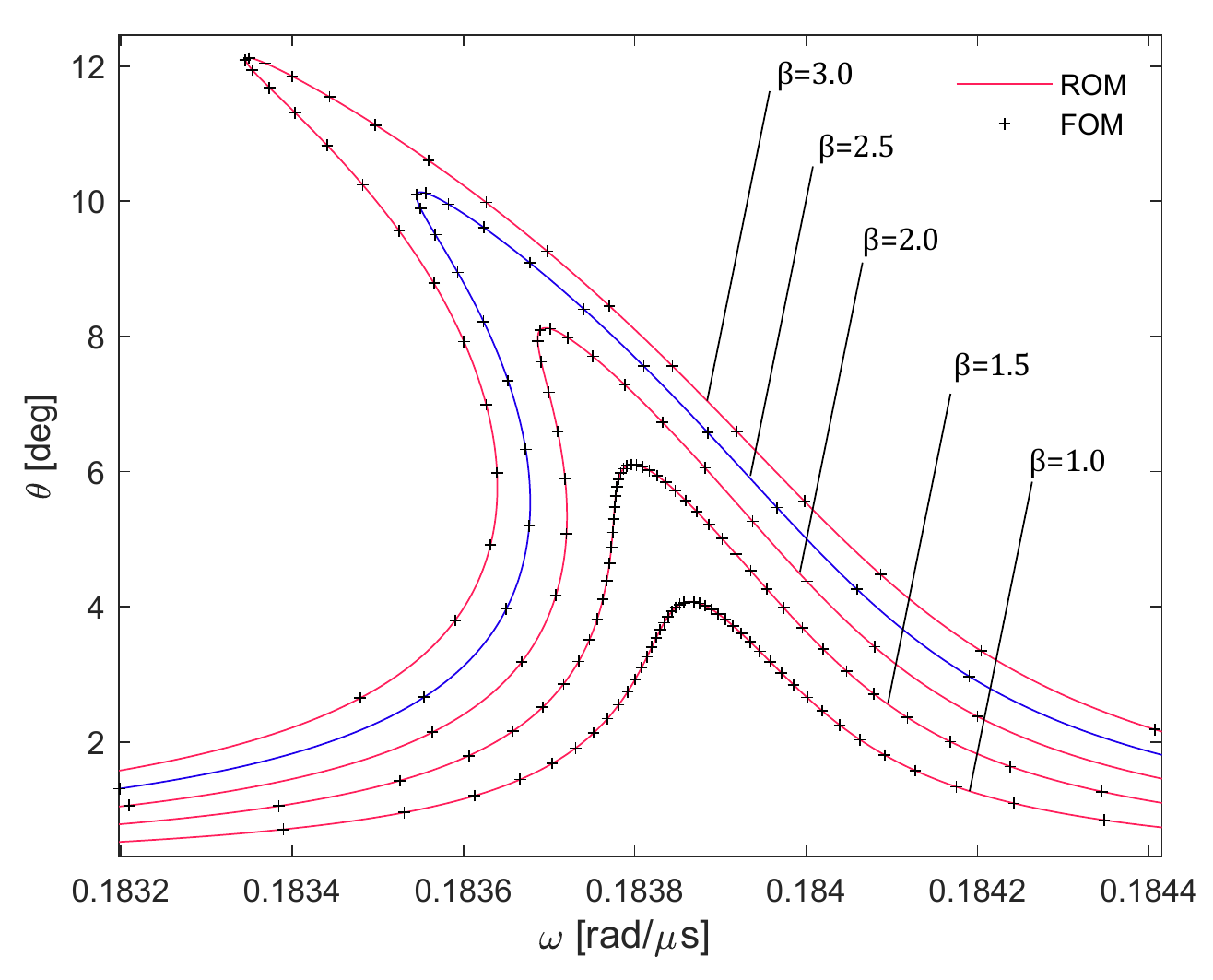}
\caption{Micromirror 2: FRFs computed with the ROMs and 8 POMs compared with the FOM solution. 
The red continuous lines denote the test ROM solutions, while 
the blue one refers to the training data. The cross markers represent the FOM solutions}
\label{fig:perseusFRF}
\end{figure}

A good balance between subspace dimension and accuracy is given by the  subspace spanned by 10 POMs
on which we perform a more extensive testing stage varying, as usual, the force multiplier 
$\beta=1, 1.5, 2, 2.5, 3\,\mu$N. 
The results plotted in Figure~\ref{fig:perseusFRF} show again the highly predictive 
capability of the POD-Galerkin ROM.

\subsection{Shallow arch with internal resonance}
\label{sec:IR}

In recent years several occurrences of complex nonlinear phenomena 
have been documented experimentally in MEMS, mainly due to their
large quality factors $Q$. 
Internal resonances (IRs) play an important role in triggering more complex motions
and facilitate energy transfer between modes. Often IRs
are strongly linked to the stability of the associated
periodic response and quasi-periodic regimes might arise as a consequence 
of Neimark-Sacker (NS) bifurcations \cite{mecc21}.
The numerical prediction of such phenomena requires an accurate stability analysis
which cannot be performed at a reasonable cost using FOMs,
while can be much more conveniently run on small ROMs using dedicated continuation tools,
as discussed in Section~\ref{sec:solROM}. 

\begin{figure}[h!]
\centering
\includegraphics[width=0.7\columnwidth]{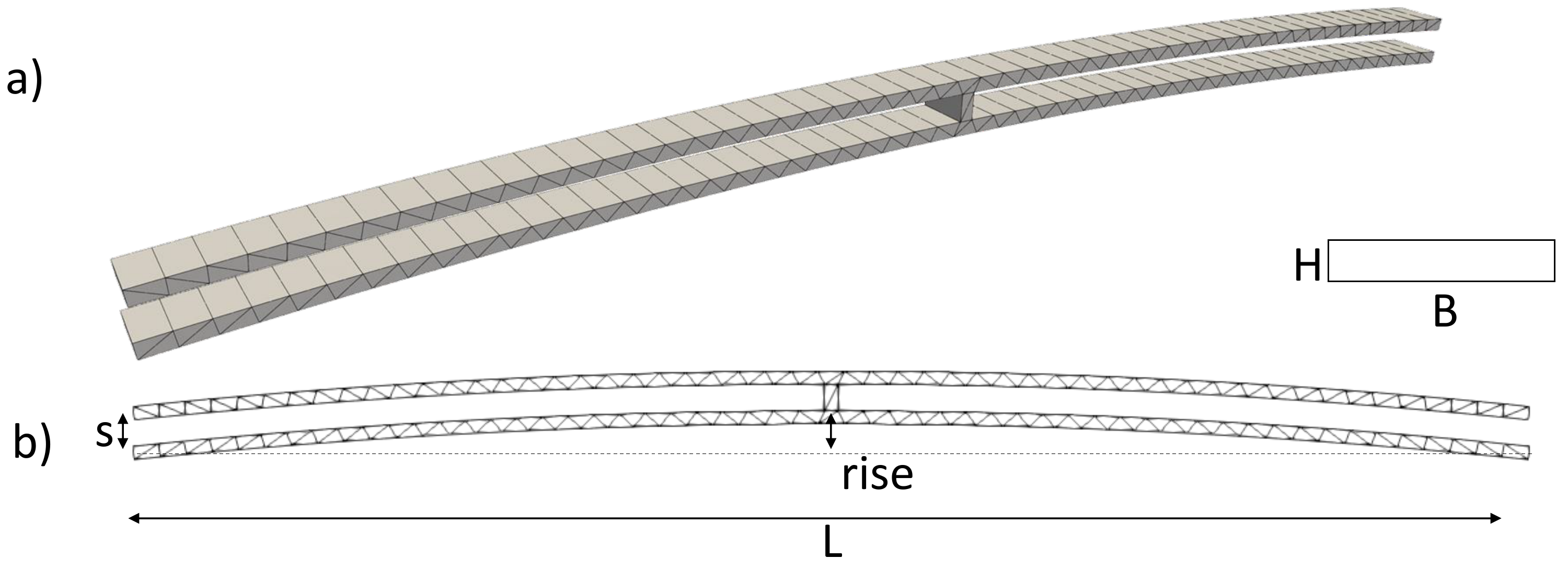}
\caption{Shallow arch. Figure a): 3D view of the FEM model. Figure b): front view and main dimensions. $B=20\,\mu$m, $H=5\,\mu$m, $L=530\,\mu$m, rise=13.4\,$\mu$m, $s=10\,\mu$m.}
\label{fig:arch1_geo}       
\end{figure}

For these reasons we include among our benchmarks a shallow
double-arch with a constant radius of curvature.
The layout, inspired by the one proposed for a bistable structure in \cite{frangi2015}, 
has been suitably designed so as to trigger a 1:2 IR.
The arch geometry and mesh are illustrated in Figure~\ref{fig:arch1_geo}. 
The mesh considered consists of quadratic wedge elements and contains 1971 nodes.

The device is made of polycrystalline silicon with density $\rho=2330$\,kg/m$^3$ and
a linear elastic Saint-Venant Kirchhoff constitutive model is assumed, with Young modulus 
$E=167000$\,MPa and Poisson coefficient $\nu=0.22$ \cite{sharpe1997measurements}. 
The first six eigenfrequencies of the modelled structure are reported in Table~\ref{tab:arch1_freq}.

\begin{table}[h!]
\centering
\begin{tabular}{ |c|c|c|c|c|c|c|} 
 \hline
Eigenmode & 1 & 2 & 3 & 4 & 5 & 6\\ 
 \hline
Frequency  [kHz]& 434.16 & 525.97 & 603.91 & 667.59 & 756.95 & 863.67 \\ 
 \hline
\end{tabular}
 \caption{First six eigenfrequencies of the MEMS arch}
\label{tab:arch1_freq}
\end{table}

\begin{figure}[h!]
\centering
\includegraphics[width = .85\linewidth]{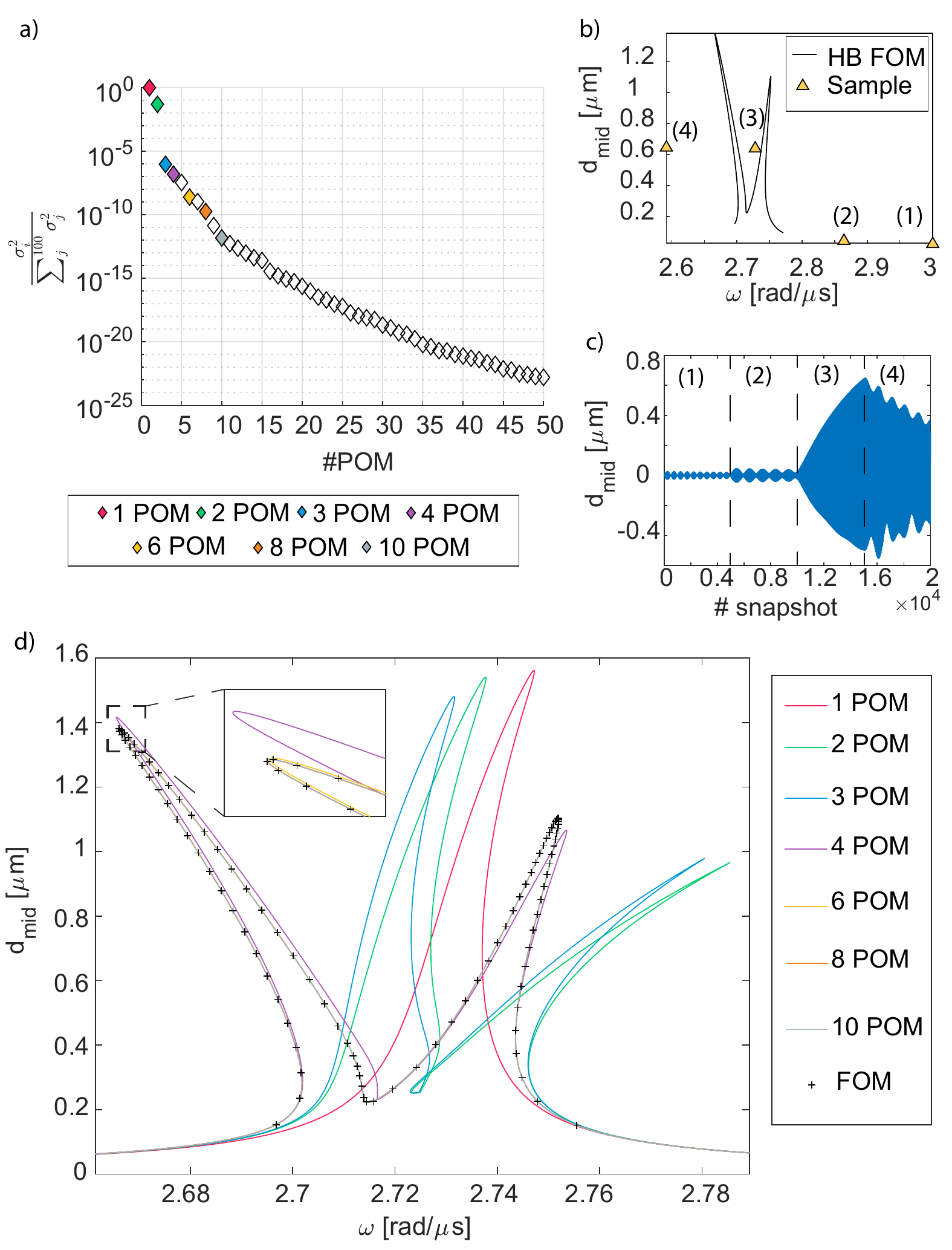}
\caption{Convergence of the hyper ROM on the MEMS arch.  Figure a): energy content of each POM. Figure b): frequency points sampled with time marching methods on the FEM model compared with the FOM FRF achieved with HB method. The solution are far from the Steady State regime and thus depart from the HB solution. Figure c): Time history of each frequency value.  Figure d: FRFs resulting from each ROM built with a increasing number of POMs. 
	The POD is tested on the training set and on a second higher forcing level.}
\label{fig:archIR_POM}
\end{figure}

The quality factor has been set to $Q=500$ and
the actuation is provided by a body force proportional to the first eigenmode 
$\bfF(t)=\bfM  \bfphi_1 \beta\cos(\omega t)$ with $\beta$ load multiplier.

In order to further stress the versatility of the POD-Galerkin ROM, 
we opt for time marching methods to simulate the FOM 
and generate the training dataset. Moreover, snapshots have been collected 
during the transient phase, far from steady-state conditions, according to the following
strategy. 
The forcing level has been fixed to $\beta=0.2\,\mu$N
and a sequence of four frequencies have been analyzed with a Newmark implicit solver
in a downward sweep. 
For each frequency, 100 cycles are simulated and a total of 20000 snapshots are collected.
Initial conditions for the global analysis are homogeneous
and the final state computed for each frequency
yields the initial conditions for the next one.
It is worth stressing that, given the $Q$
value at hand, after 100 cycles the system is still fully in a transient phase.
The frequency points and the time series of the mid-span deflection of the arch are collected 
in Figures~\ref{fig:archIR_POM}b and \ref{fig:archIR_POM}c, respectively. 
The maximum amplitudes, denoted by yellow triangles in Figure~\ref{fig:archIR_POM}b,
are indeed quite far from the steady state solutions predicted by the HB FEM.

\begin{figure}[h!]
\centering
\includegraphics[width = .7\linewidth]{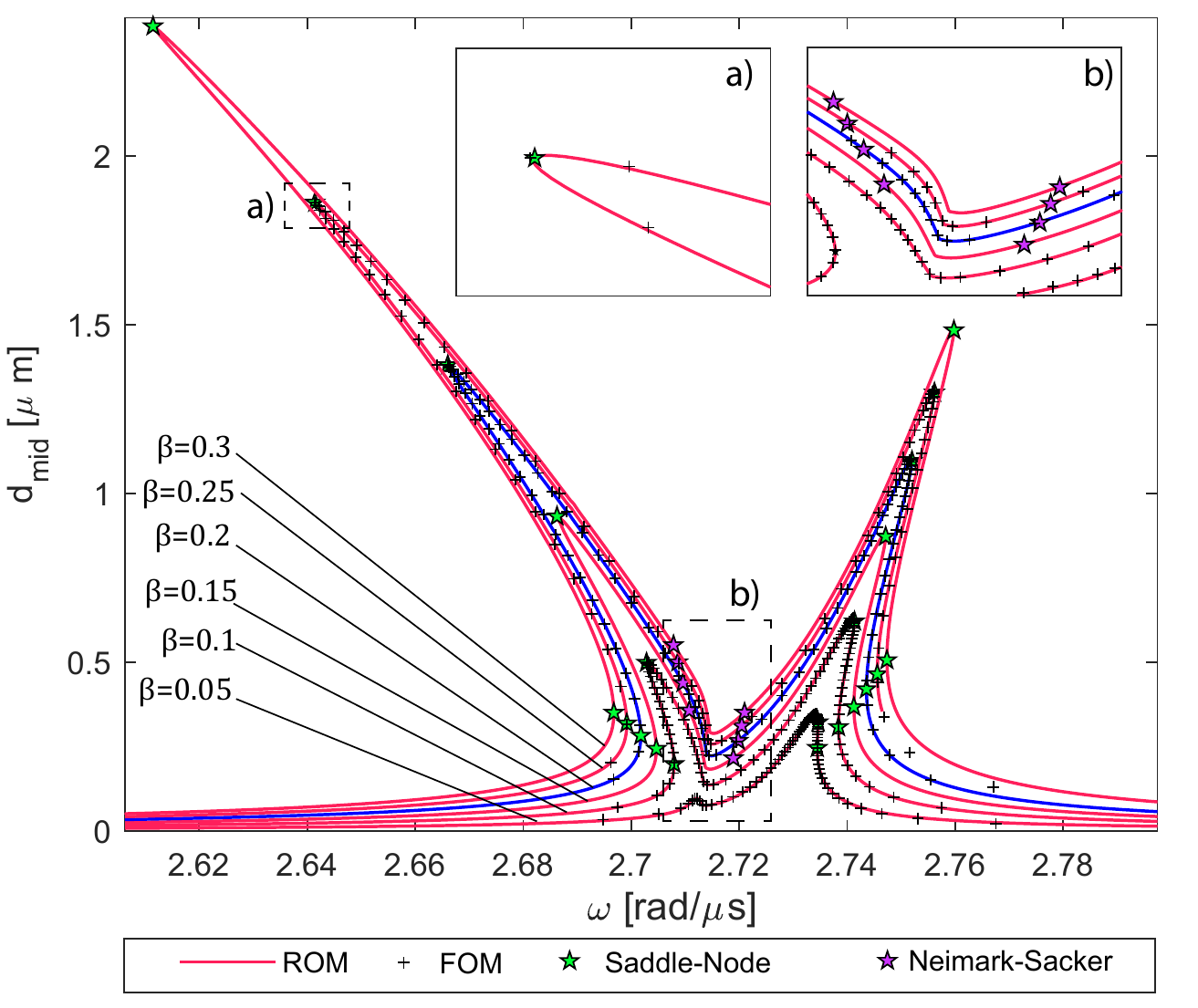}
\caption{Shallow arch: FRFs computed with the ROM (8 POMs) compared with the FOM solutions. 
The red continuous lines are ROM solutions computed in the test phase, 
the blue one marking the training curve. 
The crosses denote the FOM solutions, 
while the star markers indicate the bifurcation points. 
The green and purple stars stand for Saddle-Node and Neimark-Sacker bifurcation points, respectively}
\label{fig:archIR}
\end{figure}

The SVD computed on the snapshot matrix yields the energy distribution of Figure~\ref{fig:archIR_POM}a. 
Also in this case, although it appears that the energy is almost totally focused 
in the first two POMs,
the convergence analysis presented in Figure~\ref{fig:archIR_POM}d shows that at least 
4 POMs are required and a good convergence is achieved starting with 6 POMs.
In the following, we will consider a trial space collecting the first 8 POMs.

The chosen ROM is now tested considering different forcing levels and 
the corresponding FRFs are plotted in Figure~\ref{fig:archIR}, together with the HB-FOM solutions.
The model correctly reproduces the complex pattern of the 1:2 IR, 
as demonstrated by the shape of the frequency response displaying the two characteristic peaks.
As recalled, a key feature of the ROM is the possibility to apply the bifurcation analysis tools
discussed in Section~\ref{sec:solROM} which yield the results
of Figure~\ref{fig:archIR}. Two different classes of bifurcation points
can be identified: 
saddle-node bifurcations, that split the FRF between unstable and stable 
branches, and Neimark-Sacker bifurcations that separate stable periodic and quasi-periodic regions. 
Quasi-periodicity is a dynamic condition where the external excitation frequency of the system is paired with an incommensurate smaller frequency that modulates the amplitude of the response (see \cite{mecc21} for further details).
For a given FRF, in the region within the two Neimark-Sacker bifurcations only 
quasi-periodic, and eventually almost
chaotic, solutions \cite{mecc21} are physically meaningful.
This challenging benchmark shows that an excellent quality of the ROM can be achieved even 
with a fast training phase based on fully transient data.

\section{Electromechanical coupled problems}
\label{sec:electro}

While previous benchmarks have addressed purely mechanical problems
with geome\-trical nonlinearities, the interest in generating an optimal linear trial space
goes beyond these applications.
In MEMS applications, multiple sources of additional nonlinearities 
come from the actuation mechanism and the most typical example 
is provided by electrostatic forces that depend 
in an intrinsically nonlinear manner on the displacement field.

In most contributions addressing both geometric and electrostatic 
effects (see e.g.~\cite{younis03,younis03b} for clamped-clamped beams) 
analytical approaches or simplified structural theories are utilized,
although their application to real MEMS often leads to results that
are only qualitatively correct or need careful device-dependent calibration.
General numerical approaches are needed as MEMS might have complicated features 
that can be hardly reduced to simple models. 
However, a coupled electromechanical FOM able to simulate complex 3D structures 
is not a standard tool even for the most advanced commercial codes and
generating snapshots of the FOM solution often comes with a computational cost that may be unsuitable for practical applications.

Most importantly, a major difficulty of the POD is associated with the evaluation
of the vector of nonlinear nodal electrostatic forces (EF).
Popular data-driven algorithms 
like the Discrete Empirical Interpolation Method (DEIM) \cite{ref41,brunton} provide
an optimal reconstruction of the full nonlinear vector starting from a collection of
snapshots of the nonlinear forces.
However, the DEIM is based on the assumption that few selected entries of the 
vector can be computed at a low cost independently of the others, while in electromechanical problems, the generation of the nonlinear vector of nodal forces
has only a marginal cost with respect to the solution phase
of the electrical sub-problem, be it solved with iterative integral equation approaches
or with FE techniques.
The development of a fast algorithm to circumvent this obstacle is still an open issue.

On the contrary, a simplified way to account for EFs through POD-based models 
can be done by exploiting the same approach successfully applied with the 
implicit condensation method in \cite{jmems20reso}.

We assume that the trial space defined from snapshots given by the mechanical simu\-lations is sufficiently rich and able to represent the displacement field also for the fully coupled problem. This assumption is reasonable when 
the perturbation of the invariant manifold induced by the electrostatic couplings is moderate
and
will in general put an upper bound to the admissible EFs, i.e.\ on the bias 
voltages imposed on the electrodes.
Since the electrostatic problem is quasi-static, the EFs
depend only on the instantaneous values of the $Q_i$. 
This implies that, in a training stage, the $\bfF$ vector in eq.\eqref{eq:PPV_d2} 
due to EFs can be computed with the FOM and 
projected on the POD subspace to generate $\bfF\aPOD$ in eq.\eqref{eq:PPV_POD}
for any given combination of the weights $Q_i$ of the POMs.
The manifold of the electrostatic forces is thus pre-computed only at discrete points
in preselected admissible ranges
and is later interpolated between knots when queried during the integration of the ROM.
Moreover, many POMs (like the ``axial'' POM in Figure~\ref{fig:beam_manifold}
or high frequency bending POMs) have negligible effects on the
EFs and can be disregarded so that EFs will depend on $p_e$ electrically active 
POMs, with $p_e\ll p$ typically.

\begin{figure}[h!]
\centering
\includegraphics[width = .8\linewidth]{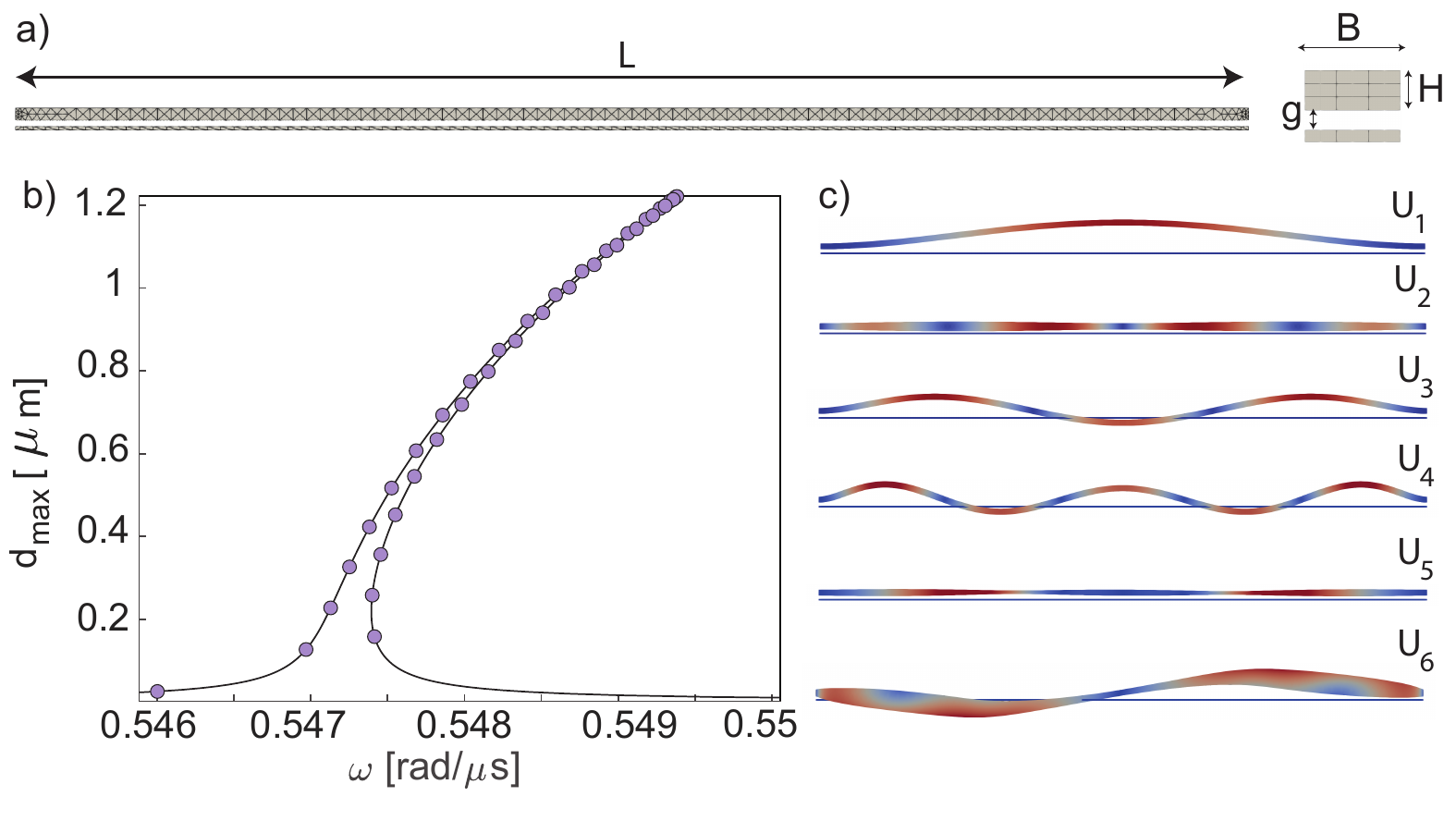}
\caption{Electromechanical problem. Figure a) Clamped clamped beam used for electromechanical simulation side and front view. Figure b) FRF used as reference for the reduction. The frequency samples used to generate the POMs are marked with circles. Figure c) fist six POMs given by the FOM snapshots.}
\label{fig:CCelettro}
\end{figure}

As a benchmark problem,
we focus on a clamped-clamped beam meshed with quadratic elements and a total of 10920 nodes. 
The dimensions of the beam are $L=1000$\,\micr, $H=10$\,\micr, $B=24$\,\micr. An electrode is
placed in front of the beam with a gap of  $g=5$\,\micr
and the voltage bias $V_{DC}+V_{AC}\cos\omega t$ is imposed between the electrode and
the beam (see Figure \ref{fig:CCelettro}a).
The quality factor has been set to $Q=10929$.
The FOM utilized for the electrostatic problem 
is a Boundary Element Code (BEM) based on integral equations accelerated with Fast Multipoles
and a total of 71844 unknowns. 
The code resorts to an iterative solver and its use in a fully coupled solution 
would have a prohibitive cost.

The ROM is trained with mechanical HB-FEM simulations setting Q=10929 and applying a body load proportional to the first eigenmode $\bfF(t)=\bfM  \bfphi_1 \beta\cos(\omega t)$ with $\beta=0.0005$. 
A total of 1850 snapshots have been collected on 37 sample frequencies
identified by circles on the FOM solution of Figure \ref{fig:CCelettro}b. 
The first 6 mechanical POMs, used to build the ROM, are depicted in Figure \ref{fig:CCelettro}c.

Consistently with the previous assumptions we can consider that only the first POM will contribute 
significantly to EFs (i.e.\ $p_e=1$). Thus, a series of electrostatic analyses are run imposing displacement fields proportional to the first POM $\mathbf{D} \approx \mathbf{U}_1 Q_1$ covering a range of 1.1\micr 
for the midspan displacement over a gap of 5\micr.
The EFs are projected on the POMs yielding the equivalent forces
expressed in $\mu$m. 
These forces are scaled by the applied potentials $V_{DC}$ and $V_{AC}$ and are modelled with a cubic polynomial as:
\begin{multline}
\label{eq:elettro}
	F\aPOD_i(Q_1,V_{DC},V_{AC},\omega,t) = \\
	=\left( V_{DC}^2 \epsilon_0 + 2 V_{DC} V_{AC}\epsilon_0 \cos(\omega t)\right) \left( \alpha^{(i)}_0+ \alpha^{(i)}_1 Q_1+ \alpha^{(i)}_2 Q_1^2  +\alpha^{(i)}_3 Q_1^3 \right)
\end{multline}
where $\alpha^{(i)}_j$ are coefficients of order $j$ associated to the force projected 
on the $i$-th POM and $\epsilon_0$ is the vacuum permittivity.
In the example considered we neglect the components proportional to $V_{AC}^2$ (since 
typically $V_{AC}\ll V_{DC}$, see e.g.\ \cite{jmems20reso,gyro2021}). 
The coefficients for the first 6 POMs are collected in Table \ref{tab:coeff_elec}.

\begin{table}[h!]
	\centering
\begin{tabular}{|c|c|c|c|c|c|c|}
	\hline
	POM & 1  & 2  & 3  &  4 & 5 & 6 \\
	\hline
	$\alpha^{(i)}_0$ & 6.8638 & -0.9609  & 3.2838 & 1.9965 & 0.2274 & 2.5945 \\
	\hline
	$\alpha^{(i)}_1$ & 0.0469 & -0.0039  & $2.86\cdot10^{-5}$ & $-3.31\cdot10^{-5}$ & 0.0021  &  0.099\\
	\hline
	$\alpha^{(i)}_2$ & $2\cdot10^{-4}$ &  $-1\cdot10^{-5}$ & $-6\cdot10^{-5}$ & $1\cdot10^{-5}$ & $1\cdot10^{-5}$ &  $3\cdot10^{-4}$ \\
	\hline
	$\alpha^{(i)}_3$ & $1\cdot10^{-6}$ & $-7\cdot10^{-8}$  & $-5\cdot10^{-7}$ & $1\cdot10^{-7}$ & $5\cdot10^{-8}$ & $1\cdot10^{-7}$  \\
	\hline
\end{tabular}
	\caption{Electromechanical problem: coefficients of the polynomial modeling	
	eq.\eqref{eq:elettro}}
	\label{tab:coeff_elec}
\end{table}

In order to provide a validation of the ROM proposed we resort 
to the commercial software MEMS+ \cite{MEMSplus}, 
which can perform electromechanical coupled simulations with structural elements. 
The model is built with 4 Timoshenko nonlinear beams with 4 nodes and the electrostatics is modelled with the method of conformal mappings. This simulation approach, though not general, is expected to be accurate for the simple MEMS tested.

\begin{table}[ht]
	\centering
	\begin{tabular}{ |c|c|c|c|c|c| } 
		\hline
		Eigenmode & 1 & 2 & 3 & 4 & 5\\ 
		\hline
		FEM model  [kHz]& 87.087 & 208.244 & 239.880 & 469.792  & 571.281 \\ 
		\hline
		MEMS+  [kHz]& 86.971 & 208.085 & 239.755 & 472.008  & 571.302 \\ 
		\hline
	\end{tabular}
	\caption{Electromechanical problem: comparison between the eigenfrequencies given by the FEM model and MEMS+}
	\label{tab:memsplus_FEM}
\end{table}

\begin{figure}[h!]
\centering
\includegraphics[width = .85\linewidth]{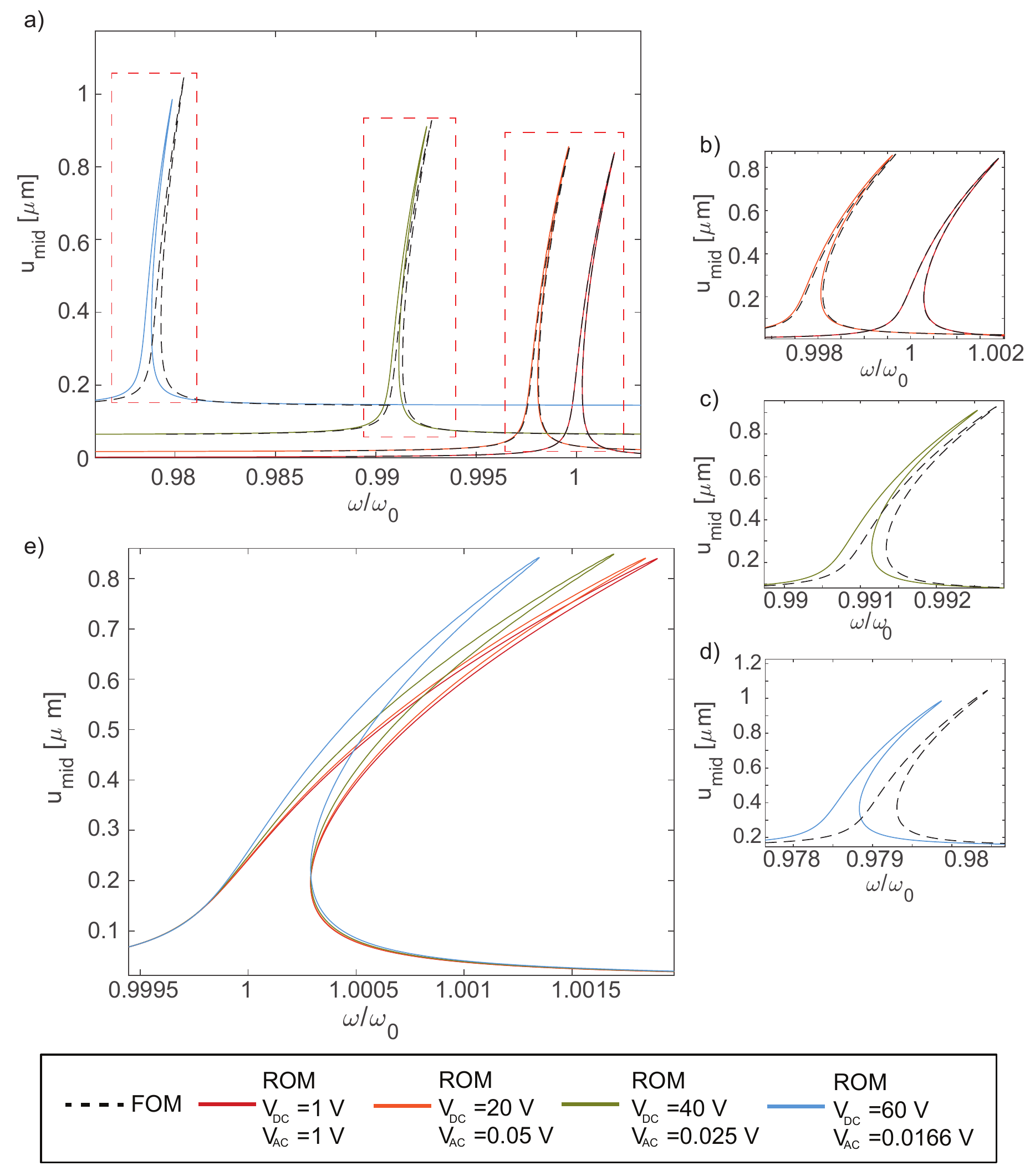}
\caption{Electromechanical problem. Figure a): FRFs computed at four different voltage bias and comparison with the MEMS+ results. Figure b) c) d): enlarged views of the FRFs computed with the ROM. Figure e): FRFs superposed filtering the shift and showing the mitigation of the hardening effect}
\label{fig:CCEcoupled}
\end{figure}

We performed several simulations considering different combinations of $V_{DC}$ and $V_{AC}$ 
giving similar peak amplitudes: $V_{DC}=1V$, $V_{AC}=1V$; $V_{DC}=20V$, $V_{AC}=0.05V$; 
$V_{DC}=40V$, $V_{AC}=0.025V$; $V_{DC}=60V$ and $V_{AC}=0.0166V$.

The results are summarized in Figure~\ref{fig:CCEcoupled}a collecting plots of the midspan deflection versus the actuation frequency.
The curves are normalized with respect to the mechanical eigenfrequency of the lowest eigenmode. The FRFs display a shift of the resonant frequency towards the left induced by the expected
electrostatic negative stiffness effect. The enlarged views in Figure~\ref{fig:CCEcoupled}b, 
\ref{fig:CCEcoupled}c and \ref{fig:CCEcoupled}d focus on specific frequency ranges for an improved 
comparison between the two classes of results.
Another relevant effect induced by EFs is given by a softening contribution 
that mitigates the hardening of the response for increasing $V_{DC}$.
To better highlight this latter effect, 
in Figure \ref{fig:CCEcoupled}e we superpose the various FRFs by filtering the mentioned shift and the static component of the displacement.

A very good agreement is achieved between the two simulation approaches
concerning the shift and the overall nonlinearity evolution. However, as expected,
for the largest values of $V_{DC}$ a quantitative mismatch  appears especially on the peak values.
Indeed, as previously remarked, in the presence of large EFs the mechanical POMs 
might not guarantee an accurate representation of mechanical displacements. The fact that 
this disagreement occurs only at very large voltages unusual in MEMS application
is however a strong validation of the proposed approach.

\section{Conclusions}
\label{sec:conclusions}

We have investigated various applications of POD-Galerkin ROMs for the sake of efficiently simulating MEMS devices. 
In parti\-cular, we have demonstrated how POD, despite being a {\it linear} reduction
technique, can tackle very efficiently and accurately highly nonlinear features common in MEMS applications: large rotations, and large geometrical transformations in general, 
internal resonances and electrostatic nonlinearities.
Geometrical nonlinearities, leading to polynomial terms up to cubic order, have been 
reduced through an exact projection onto the subspace spanned by the POMs, 
while electrostatics has been modeled resorting to precomputed manifolds in terms of the amplitudes of the electrically active POMs. 

We have tested extensively the reliability of the assumed linear trial space in challenging applications focusing on resonators, micromirrors and arches also displaying internal resonances.
We have discussed several options to generate the matrix of snapshots 
using both classical time marching schemes
and more advanced Harmonic Balance approaches. It has been shown that the method is robust
and that a POD-Galerkin ROM can be trained also with transient time simulations far from steady state conditions.
This might indeed be the only viable option in many applications, considering the cost and complexity 
of HB methods and the large quality factors typical of MEMS 
that prevent reaching steady state conditions with time-marching schemes.

In order to provide a deeper insight into the POD approach, 
we have also shown the similarity of the POMs extracted 
with the reconstruction vectors of the DP approach that relies on the invariant manifold theory
and verified that the trajectories predicted by the POD-Galerkin ROM perfectly lie on the manifolds of the DP.
This result further strengthens the interpretation of the POMs as the best linear approximation of the Nonlinear Normal Modes in the least square sense \cite{amabili07}.

Another relevant feature of a POD-Galerkin ROM is that is allows to apply continuation
methods and stability analysis of the dynamic solution, usually infeasible 
in FOM analyses.
This provides the possibility to compute directly periodic solutions, trace the full FRF
with both stable and unstable branches, and locate bifurcation points
by resorting to continuation codes available in the literature.

The whole set of challenging benchmarks developed seems to suggest
a high potentiality, possibly so far underestimated,  
of the POD for this specific class of applications and 
stresses the reliability of the technique and its strong predictive ability.


\bibliography{biblio}

\begin{thebibliography}{10}
\expandafter\ifx\csname url\endcsname\relax
  \def\url#1{\texttt{#1}}\fi
\expandafter\ifx\csname urlprefix\endcsname\relax\def\urlprefix{URL }\fi
\expandafter\ifx\csname href\endcsname\relax
  \def\href#1#2{#2} \def\path#1{#1}\fi

\bibitem{rega2005}
G.~Rega, H.~Troger, Dimension reduction of dynamical systems: methods, models,
  applications, Nonlinear Dynamics 41~(1) (2005) 1--15.

\bibitem{kerschen2005}
G.~Kerschen, J.-c. Golinval, A.~F. Vakakis, L.~A. Bergman, The method of proper
  orthogonal decomposition for dynamical characterization and order reduction
  of mechanical systems: an overview, Nonlinear Dynamics 41~(1) (2005)
  147--169.

\bibitem{gordon2008}
J.~J. Hollkamp, R.~W. Gordon, Reduced-order models for nonlinear response
  prediction: Implicit condensation and expansion, Journal of Sound and
  Vibration 318~(4-5) (2008) 1139--1153.

\bibitem{mignolet}
M.~P. Mignolet, A.~Przekop, S.~A. Rizzi, S.~M. Spottswood, A review of
  indirect/non-intrusive reduced order modeling of nonlinear geometric
  structures, Journal of Sound and Vibration 332~(10) (2013) 2437--2460.

\bibitem{ijnm19}
A.~Frangi, G.~Gobat, Reduced order modelling of the non-linear stiffness in
  mems resonators, International Journal of Non-Linear Mechanics 116 (2019)
  211--218.

\bibitem{jmems20reso}
V.~Zega, G.~Gattere, S.~Koppaka, A.~Alter, G.~D. Vukasin, A.~Frangi, T.~W.
  Kenny, Numerical modelling of non-linearities in mems resonators, Journal of
  Microelectromechanical Systems 29~(6) (2020) 1443--1454.

\bibitem{nayfeh95}
A.~H. Nayfeh, D.~T. Mook, P.~Holmes, Nonlinear oscillations (1980).

\bibitem{nayfeh89}
A.~H. Nayfeh, B.~Balachandran, Modal interactions in dynamical and structural
  systems (1989).

\bibitem{krylov2011}
S.~Krylov, B.~R. Ilic, S.~Lulinsky, Bistability of curved microbeams actuated
  by fringing electrostatic fields, Nonlinear Dynamics 66~(3) (2011) 403--426.

\bibitem{Czaplewski2019}
D.~A. Czaplewski, S.~Strachan, O.~Shoshani, S.~W. Shaw, D.~L{\'o}pez,
  Bifurcation diagram and dynamic response of a mems resonator with a 1: 3
  internal resonance, Applied Physics Letters 114~(25) (2019) 254104.

\bibitem{Houri2019}
S.~Houri, D.~Hatanaka, M.~Asano, R.~Ohta, H.~Yamaguchi, Limit cycles and
  bifurcations in a nonlinear mems resonator with a 1: 3 internal resonance,
  Applied Physics Letters 114~(10) (2019) 103103.

\bibitem{avoort10}
C.~Van~der Avoort, R.~Van~der Hout, J.~Bontemps, P.~Steeneken, K.~Le~Phan,
  R.~Fey, J.~Hulshof, J.~Van~Beek, Amplitude saturation of mems resonators
  explained by autoparametric resonance, Journal of Micromechanics and
  Microengineering 20~(10) (2010) 105012.

\bibitem{ruzziconi2021two}
L.~Ruzziconi, N.~Jaber, L.~Kosuru, M.~L. Bellaredj, M.~I. Younis, Two-to-one
  internal resonance in the higher-order modes of a mems beam: Experimental
  investigation and theoretical analysis via local stability theory,
  International Journal of Non-Linear Mechanics 129 (2021) 103664.

\bibitem{nitzan16}
S.~H. Nitzan, P.~Taheri-Tehrani, M.~Defoort, S.~Sonmezoglu, D.~A. Horsley,
  Countering the effects of nonlinearity in rate-integrating gyroscopes, IEEE
  Sensors Journal 16~(10) (2016) 3556--3563.

\bibitem{seshia17}
A.~Ganesan, C.~Do, A.~Seshia, Phononic frequency comb via intrinsic three-wave
  mixing, Physical Review Letters 118~(3) (2017) 033903.

\bibitem{seshia18}
A.~Ganesan, C.~Do, A.~Seshia, Phononic frequency comb via three-mode parametric
  resonance, Applied Physics Letters 112~(2) (2018) 021906.

\bibitem{harvester12}
L.~Xiong, L.~Tang, B.~R. Mace, Internal resonance with commensurability induced
  by an auxiliary oscillator for broadband energy harvesting, Applied Physics
  Letters 108~(20) (2016) 203901.

\bibitem{KerschenHB2015}
T.~Detroux, L.~Renson, L.~Masset, G.~Kerschen, The harmonic balance method for
  bifurcation analysis of large-scale nonlinear mechanical systems, Computer
  Methods in Applied Mechanics and Engineering 296 (2015) 18--38.

\bibitem{actuators21}
A.~Opreni, N.~Boni, R.~Carminati, A.~Frangi, Analysis of the nonlinear response
  of piezo-micromirrors with the harmonic balance method, in: Actuators,
  Vol.~10, MDPI, 2021, p.~21.

\bibitem{rizzi2003}
A.~A. Muravyov, S.~A. Rizzi, Determination of nonlinear stiffness with
  application to random vibration of geometrically nonlinear structures,
  Computers \& Structures 81~(15) (2003) 1513--1523.

\bibitem{givois}
A.~Givois, A.~Grolet, O.~Thomas, J.-F. De{\"u}, On the frequency response
  computation of geometrically nonlinear flat structures using reduced-order
  finite element models, Nonlinear Dynamics 97~(2) (2019) 1747--1781.

\bibitem{cyrilstep}
A.~Vizzaccaro, A.~Givois, P.~Longobardi, Y.~Shen, J.-F. De{\"u}, L.~Salles,
  C.~Touz{\'e}, O.~Thomas, Non-intrusive reduced order modelling for the
  dynamics of geometrically nonlinear flat structures using three-dimensional
  finite elements, Computational Mechanics 66~(6) (2020) 1293--1319.

\bibitem{gyro2021}
G.~Gobat, V.~Zega, P.~Fedeli, L.~Guerinoni, C.~Touz\'e, A.~Frangi, Reduced
  order modelling and experimental validation of a mems gyroscope
  test-structure exhibiting 1:2 internal resonance, Scientific Reports (2021 -
  in-press).

\bibitem{nicolaidou2020indirect}
E.~Nicolaidou, T.~L. Hill, S.~A. Neild, Indirect reduced-order modelling: using
  nonlinear manifolds to conserve kinetic energy, Proceedings of the Royal
  Society A 476~(2243) (2020) 20200589.

\bibitem{alfiobook}
A.~Quarteroni, A.~Manzoni, F.~Negri, Reduced basis methods for partial
  differential equations: an introduction, Vol.~92, Springer, 2015.

\bibitem{lu2019}
K.~Lu, Y.~Jin, Y.~Chen, Y.~Yang, L.~Hou, Z.~Zhang, Z.~Li, C.~Fu, Review for
  order reduction based on proper orthogonal decomposition and outlooks of
  applications in mechanical systems, Mechanical Systems and Signal Processing
  123 (2019) 264--297.

\bibitem{md2016}
O.~Weeger, U.~Wever, B.~Simeon, On the use of modal derivatives for nonlinear
  model order reduction, International Journal for Numerical Methods in
  Engineering 108~(13) (2016) 1579--1602.

\bibitem{QM2017}
S.~Jain, P.~Tiso, J.~B. Rutzmoser, D.~J. Rixen, A quadratic manifold for model
  order reduction of nonlinear structural dynamics, Computers \& Structures 188
  (2017) 80--94.

\bibitem{QMgen2017}
J.~B. Rutzmoser, D.~J. Rixen, P.~Tiso, S.~Jain, Generalization of quadratic
  manifolds for reduced order modeling of nonlinear structural dynamics,
  Computers \& Structures 192 (2017) 196--209.

\bibitem{mahdiabadi2021non}
M.~K. Mahdiabadi, P.~Tiso, A.~Brandt, D.~J. Rixen, A non-intrusive model-order
  reduction of geometrically nonlinear structural dynamics using modal
  derivatives, Mechanical Systems and Signal Processing 147 (2021) 107126.

\bibitem{rosenberg}
R.~M. Rosenberg, The normal modes of nonlinear n-degree-of-freedom systems
  (1962).

\bibitem{shawNNM}
S.~Shaw, An invariant manifold approach to nonlinear normal modes of
  oscillation, Journal of Nonlinear Science 4~(1) (1994) 419--448.

\bibitem{touze2004}
C.~Touz{\'e}, O.~Thomas, A.~Chaigne, Hardening/softening behaviour in
  non-linear oscillations of structural systems using non-linear normal modes,
  Journal of Sound and Vibration 273~(1-2) (2004) 77--101.

\bibitem{touze2006}
C.~Touz{\'e}, M.~Amabili, Nonlinear normal modes for damped geometrically
  nonlinear systems: Application to reduced-order modelling of harmonically
  forced structures, Journal of Sound and Vibration 298~(4-5) (2006) 958--981.

\bibitem{touze}
C.~Touz{\'e}, Normal form theory and nonlinear normal modes: theoretical
  settings and applications, in: Modal analysis of nonlinear mechanical
  systems, Springer, 2014, pp. 75--160.

\bibitem{haller16}
G.~Haller, S.~Ponsioen, Nonlinear normal modes and spectral submanifolds:
  existence, uniqueness and use in model reduction, Nonlinear Dynamics 86~(3)
  (2016) 1493--1534.

\bibitem{haller18}
S.~Ponsioen, T.~Pedergnana, G.~Haller, Automated computation of autonomous
  spectral submanifolds for nonlinear modal analysis, Journal of Sound and
  Vibration 420 (2018) 269--295.

\bibitem{KerschenNNM}
L.~Renson, G.~Kerschen, B.~Cochelin, Numerical computation of nonlinear normal
  modes in mechanical engineering, Journal of Sound and Vibration 364 (2016)
  177--206.

\bibitem{vizza2020}
A.~Vizzaccaro, Y.~Shen, L.~Salles, J.~Blaho{\v{s}}, C.~Touz{\'e}, Direct
  computation of nonlinear mapping via normal form for reduced-order models of
  finite element nonlinear structures, arXiv preprint arXiv:2009.12145 (2020).

\bibitem{nld21}
A.~Opreni, A.~Vizzaccaro, A.~Frangi, C.~Touz{\'e}, Model order reduction based
  on direct normal form: Application to large finite element mems structures
  featuring internal resonance, Nonlinear Dynamics (2021 - in press).
\newblock \href {https://doi.org/10.1007/s11071-021-06641-7}
  {\path{doi:10.1007/s11071-021-06641-7}}.

\bibitem{kerschen2002}
G.~Kerschen, J.-C. Golinval, Physical interpretation of the proper orthogonal
  modes using the singular value decomposition, Journal of Sound and Vibration
  249~(5) (2002) 849--865.

\bibitem{sampaio}
R.~Sampaio, C.~Soize, Remarks on the efficiency of pod for model reduction in
  non-linear dynamics of continuous elastic systems, International Journal for
  Numerical Methods in Engineering 72~(1) (2007) 22--45.

\bibitem{amabili07}
M.~Amabili, C.~Touz{\'e}, Reduced-order models for nonlinear vibrations of
  fluid-filled circular cylindrical shells: comparison of pod and asymptotic
  nonlinear normal modes methods, Journal of Fluids and Structures 23~(6)
  (2007) 885--903.

\bibitem{tiso}
P.~Tiso, D.~J. Rixen, Reduction methods for mems nonlinear dynamic analysis,
  in: Nonlinear Modeling and Applications, Volume 2, Springer, 2011, pp.
  53--65.

\bibitem{korvink}
J.~S. Han, E.~B. Rudnyi, J.~G. Korvink, Efficient optimization of transient
  dynamic problems in mems devices using model order reduction, Journal of
  Micromechanics and Microengineering 15~(4) (2005) 822.

\bibitem{malvern}
L.~E. Malvern, Introduction to the Mechanics of a Continuous Medium, no.
  Monograph, 1969.

\bibitem{barrault2004anempirical}
M.~Barrault, Y.~Maday, N.~C. Nguyen, A.~T. Patera, An 'empirical interpolation'
  method: Application to efficient reduced-basis discretization of partial
  differential equations, Comptes Rendus Math\'ematique de l'Acad\'emie des
  Sciences 339~(9) (2004) 667--672.

\bibitem{chaturantabut2010nonlinear}
S.~Chaturantabut, D.~C. Sorensen, Nonlinear model reduction via discrete
  empirical interpolation, SIAM Journal on Scientific Computing 32~(5) (2010)
  2737--2764.

\bibitem{maday2008ageneral}
Y.~Maday, N.~Nguyen, A.~Patera, G.~Pau, A general multipurpose interpolation
  procedure: The magic points, Communications on Pure and Applied Analysis 8
  (2008) 383--404.

\bibitem{Doedel}
E.~J. Doedel, A.~R. Champneys, F.~Dercole, T.~F. Fairgrieve, Y.~A. Kuznetsov,
  B.~Oldeman, R.~Paffenroth, B.~Sandstede, X.~Wang, C.~Zhang, Auto-07p:
  Continuation and bifurcation software for ordinary differential equations
  (2007).

\bibitem{Guillot1}
L.~Guillot, B.~Cochelin, C.~Vergez, A taylor series-based continuation method
  for solutions of dynamical systems, Nonlinear Dynamics 98~(4) (2019)
  2827--2845.

\bibitem{Guillot2}
L.~Guillot, A.~Lazarus, O.~Thomas, C.~Vergez, B.~Cochelin, A purely frequency
  based floquet-hill formulation for the efficient stability computation of
  periodic solutions of ordinary differential systems, Journal of Computational
  Physics 416 (2020) 109477.

\bibitem{NVLIB}
M.~Krack, J.~Gross, Harmonic balance for nonlinear vibration problems, Vol.~1,
  Springer, 2019.

\bibitem{COCO}
H.~Dankowicz, F.~Schilder, Recipes for continuation, SIAM, 2013.

\bibitem{BifurcationKit}
R.~Veltz, Bifurcationkit. jl, Ph.D. thesis, Inria Sophia-Antipolis (2020).

\bibitem{cochelin2009high}
B.~Cochelin, C.~Vergez, A high order purely frequency-based harmonic balance
  formulation for continuation of periodic solutions, Journal of Sound and
  Vibration 324~(1-2) (2009) 243--262.

\bibitem{buza2021using}
G.~Buza, S.~Jain, G.~Haller, Using spectral submanifolds for optimal mode
  selection in nonlinear model reduction, Proceedings of the Royal Society A
  477~(2246) (2021) 20200725.

\bibitem{haro2006parameterization1}
A.~Haro, R.~de~la Llave, A parameterization method for the computation of
  invariant tori and their whiskers in quasi-periodic maps: numerical
  algorithms, Discrete \& Continuous Dynamical Systems-B 6~(6) (2006) 1261.

\bibitem{haro2006parameterization2}
A.~Haro, R.~de~la Llave, A parameterization method for the computation of
  invariant tori and their whiskers in quasi-periodic maps: rigorous results,
  Journal of Differential Equations 228~(2) (2006) 530--579.

\bibitem{haro2016parameterization}
A.~Haro, M.~Canadell, J.-L. Figueras, A.~Luque, J.-M. Mondelo, The
  parameterization method for invariant manifolds, Applied Mathematical
  Sciences 195 (2016).

\bibitem{jain2021compute}
S.~Jain, G.~Haller, How to compute invariant manifolds and their reduced
  dynamics in high-dimensional finite-element models?, arXiv preprint
  arXiv:2103.10264 (2021).

\bibitem{jmems04}
A.~Corigliano, B.~De~Masi, A.~Frangi, C.~Comi, A.~Villa, M.~Marchi, Mechanical
  characterization of polysilicon through on-chip tensile tests, Journal of
  Microelectromechanical Systems 13~(2) (2004) 200--219.

\bibitem{ARlenses}
Laser beam scanning,
  \url{https://www.st.com/content/st_com/en/about/innovation---technology/laser-beam-scanning.html}.

\bibitem{jmemsmirror}
A.~Frangi, A.~Opreni, N.~Boni, P.~Fedeli, R.~Carminati, M.~Merli, G.~Mendicino,
  Nonlinear response of pzt-actuated resonant micromirrors, Journal of
  Microelectromechanical Systems 29~(6) (2020) 1421--1430.

\bibitem{hopcroft2010}
M.~A. Hopcroft, W.~D. Nix, T.~W. Kenny, What is the young's modulus of
  silicon?, Journal of Microelectromechanical Systems 19~(2) (2010) 229--238.

\bibitem{mecc21}
G.~Gobat, L.~Guillot, A.~Frangi, B.~Cochelin, C.~Touz{\'e}, Backbone curves,
  neimark-sacker boundaries and appearance of quasi-periodicity in nonlinear
  oscillators: application to 1: 2 internal resonance and frequency combs in
  mems, Meccanica (2021) 1--33.

\bibitem{frangi2015}
A.~Frangi, B.~De~Masi, F.~Confalonieri, S.~Zerbini, Threshold shock sensor
  based on a bistable mechanism: design, modeling, and measurements, Journal of
  Microelectromechanical Systems 24~(6) (2015) 2019--2026.

\bibitem{sharpe1997measurements}
W.~N. Sharpe, B.~Yuan, R.~Vaidyanathan, R.~L. Edwards, Measurements of young's
  modulus, poisson's ratio, and tensile strength of polysilicon, in:
  Proceedings IEEE the tenth annual international workshop on micro electro
  mechanical systems. An investigation of micro structures, sensors, actuators,
  machines and robots, IEEE, 1997, pp. 424--429.

\bibitem{younis03}
M.~I. Younis, E.~M. Abdel-Rahman, A.~Nayfeh, A reduced-order model for
  electrically actuated microbeam-based mems, Journal of Microelectromechanical
  Systems 12~(5) (2003) 672--680.

\bibitem{younis03b}
M.~I. Younis, A.~Nayfeh, A study of the nonlinear response of a resonant
  microbeam to an electric actuation, Nonlinear Dynamics 31~(1) (2003) 91--117.

\bibitem{ref41}
M.~Barrault, Y.~Maday, N.~C. Nguyen, A.~T. Patera, An ' empirical interpolation
  ' method: application to efficient reduced-basis discretization of partial
  differential equations, Comptes Rendus Mathematique 339~(9) (2004) 667--672.

\bibitem{brunton}
S.~L. Brunton, J.~N. Kutz, Data-driven science and engineering: Machine
  learning, dynamical systems, and control, Cambridge University Press, 2019.

\bibitem{MEMSplus}
J.~Nazdrowicz, A.~Napieralski, Modelling, simulations and performance analysis
  of mems vibrating gyroscope in coventor mems+ environment, in: 2019 20th
  International Conference on Thermal, Mechanical and Multi-Physics Simulation
  and Experiments in Microelectronics and Microsystems (EuroSimE), IEEE, 2019,
  pp. 1--5.

\end{thebibliography}

\appendix 

\section{Analysis of computational performances}
\label{sec:Annex}
In this Appendix we discuss the computational performances of the FOM and of the offline-online 
stages of the ROM for the mechanical examples of Section \ref{sec:mech}.
These data further stress the efficiency of the ROM and provide the trend of the achieved speed-up.

The comparison between the cost of the FOM and of the online stage of the ROM is reported in Table \ref{tab:annex_FOMon}, while an analysis of the offline stage is reported in Table \ref{tab:annex_off}. 
All the simulations have been run on a workstation with AMD Ryzen 5 1600 Six-Core Processor 3.20 GHz with 64 GB RAM.
For HB methods the computational time only includes the calculation of the harmonic components and it does not account for the reconstruction of the time history over the period, while
for time-marching methods the cost of all the steps of the time history is provided.

\begin{table}[h!]
	\centering
	\begin{tabular}{|c|c|c|c|c|c|}
		\hline
		Application &  FOM-m & $T_{\text{FOM}}$ [$10^6$ s]  & $p$&  $T_{\text{online}} [s]$   & 
		$T_{\text{FOM}}/T_{\text{online}}$  \\
		\hline
		C-C Beam & HB (9)  & 0.21 & 6 & 42 & 5023 \\
		\hline
		C-C Beam & TM-SS (50)  & 1.35 & 6 & 42 & 32143 \\
			\hline
		Micromirror 1 & HB (5) & 0.44 & 8 &452& 969\\
		\hline
		Micromirror 2& HB (7) &  0.2 & 10 &  220 & 918\\
		\hline
		Arch & HB (9) & 0.12 & 8 & 250 & 480 \\
		\hline
		Arch & TM-SS (50)      & 9 & 8 & 250 & 36000\\
		\hline
	\end{tabular}
	\caption{Cost of the FOM and of the online stage of the ROM considering 1000 frequency instances for a fixed forcing level}
	\label{tab:annex_FOMon}
\end{table}

In Table~\ref{tab:annex_FOMon} 
the FOM-m column specifies the solution technique used for the full order model
according to the conventions introduced in Section \ref{sec:mech}. 
Within brackets we provide the number of time steps in one period 
for time marching methods (TM)  or the number of harmonic components for HB methods.
$T_{\text{FOM}}$ is the average computing time required to obtain a FOM solution for 1000 frequencies and a single forcing level, which is equivalent to a finely sampled FRF with continuation methods. 
The $p$  column provides the number of POMs used in the reduction. 
$T_{\text{online}}$ is the average time required to compute the ROM solution for 1000 frequencies and a single forcing level. The last column provides the speed-up
defined in this case as the ration between $T_{\text{FOM}}$ and $T_{\text{online}}$. Indeed, when the offline cost is not an issue and one is mainly interested in the online phase, this is the 
most important performance indicator.
One comment is worth stressing concerning the FOM solved with time marching methods.
It appears that the time marching approach is not suitable for MEMS applications 
when the steady state regime is needed, as is the case herein. Indeed this is due to the 
large quality factors $Q$ typical of MEMS. 
For the C-C beam and the arch resonator we consider that the steady state is reached after $6 Q$ periods and, despite the fact that the $Q$ factors considered are unrealistically small, the computing cost explodes.

\begin{table}[h!]
	\centering
	\begin{tabular}{|c|c|c|c|c|c|c|c|}
		\hline
		\rule{0pt}{12pt}
		Application &  FOM-m&  $\#\omega$ & $m$ & $T_{\text{snap}} [s]$ & $T_{\text{offline}} [s]$ & $\frac{T_{\text{FOM}}}{T_{\text{online}}+T_{\text{offline}}}$  \\
		\hline
		C-C Beam& HB (9) & 10  &500 & 2110 & 2 301  &90\\
		\hline
		\rule{0pt}{12pt}
		C-C Beam& TM-SS  (50)  & 2  &  1242 & 720 & 770  &259\\
		\hline
		\rule{0pt}{12pt}
		C-C Beam& TM-SS  (50)  & 10  &  6210 & 3600 & 4014  &52\\
		\hline
		\rule{0pt}{12pt}
		C-C Beam& TM-SS  (50)  & 20  &  12420 & 7200 & 8206  & 25\\
		\hline
		\rule{0pt}{12pt}
		C-C Beam& TM-TR  (100) & 1  &  1000 & 90 & 137  & 1178\\
		\hline
		Micromirror 1& HB (5)  & 50  & 5000 &3026 & 25526 & 16 \\
		\hline
		Micromirror 2& HB (7) & 40 & 2000 & 208& 8538 & 23\\
		\hline
		\rule{0pt}{12pt}
		Arch & TM  (50)  & 4  & 20000& 484 & 2214 &48 \\
		\hline
	\end{tabular}
	\caption{Analysis of the offline stage}
	\label{tab:annex_off}
\end{table}

In Table~\ref{tab:annex_off} we focus on the contrary on the offline stage and provide a second 
possible measure of speed-up. 
As before, the FOM-m column specifies the solution technique utilized for the full order model.
$\# \omega$  is the number of frequency instances used in the training
and  $m$  is the total number of snapshots used in the SVD (including both frequency and forcing variations). 
The $T_{\text{snap}}$  column provides the time required to compute the snapshots with the FOM.
$T_{\text{offline}}$ is the sum of $T_{\text{snap}}$, of the time required to perform the SVD 
decomposition (package used ARPACK in \texttt{FORTRAN}) and of the cost of the projection onto the ROM subspace. 
Finally we report a second possible speed measure $T_{\text{FOM}}/(T_{\text{online}}+T_{\text{offline}})$ that includes also the impact of $T_{\text{offline}}$. The $T_{\text{FOM}}$ and the $T_{\text{online}}$ are the one in Table~\ref{tab:annex_FOMon}, HB method is considered in $T_{\text{FOM}}$. This alternative speed measure, shows that $T_{\text{offline}}$ represents a significant part of the computational effort, nevertheless the ROM is still convenient. The time spent in the offline stage may overcome the gain achieved in the online stage when few parameter queries are simulated with the ROM. However in MEMS applications, considering all the features mentioned in this work (e.g high Q factor, large models, multiple parameters to span etc.) the computation gain is always much greater than 1.

Let us consider the C-C beam case where 
we compare the offline cost of four different conditions: HB FOM, 5,10 and  20 frequency samples with snapshots taken close to steady state (SS, see Figure \ref{fig:beam_samples}a)  and one frequency with transient time series (TR, see Figure \ref{fig:beam_samples}b) generated with TM methods.
This highlights that TM methods represent an appealing alternative 
as a limited number of snapshots sampled in a fully transient state still allows 
to identify a proper subspace, as pointed out in Section \ref{sec:ccbeam}.

\end{document}